\documentclass[12pt,twoside,a4paper]{amsart}

\usepackage{amssymb}
\usepackage{latexsym}
\usepackage{amsmath}
\usepackage{euscript}
\usepackage{graphics}

\newcommand\tA{\widetilde{A}}

\newcommand\gt{{\mathfrak{t}}}
\newcommand\gH{{\mathfrak{H}}}
\newcommand\gN{{\mathfrak{N}}}
\newcommand\gR{{\mathfrak{R}}}
\newcommand\gS{{\mathfrak{S}}}

\newcommand\T{{\mathbb{T}}}

\newcommand\ff{\varphi}
\newcommand\kk{\varkappa}
\newcommand\tg{\tan}
\newcommand\ctg{\cot}
\newcommand\equ{{\Longleftrightarrow}}
\newcommand\then{\Longrightarrow}
\newcommand\eset{\varnothing}

\newcommand\cH{{\mathcal{H}}}
\newcommand\cD{{\mathcal{D}}}
\newcommand\cC{{\mathcal{C}}}

\def\Llr{\Longleftrightarrow}

\def\wt#1{{{\widetilde #1} }}

\def\bm\chi{\mbox{\boldmath$\chi$}}

\def\Ext{{\rm Ext\,}}
\def\Extp{{\rm Extp\,}}

\def\ran{{\rm ran\,}}
\def\dom{{\rm dom\,}}

\def\dim{{\rm dim\,}}

\let\xker=\ker \def\ker{{\xker\,}}

\DeclareMathOperator\re{Re}
\DeclareMathOperator\im{Im}

\newcommand\dup{\{T_1,T_2\}}

\newtheorem{theorem}{Theorem}[section]
\newtheorem{proposition}[theorem]{Proposition}
\newtheorem{corollary}[theorem]{Corollary}
\newtheorem{lemma}[theorem]{Lemma}
\newtheorem{definition}[theorem]{Definition}
\theoremstyle{definition}

\newtheorem{remark}[theorem]{Remark}

\numberwithin{equation}{section}

\begin{document}
\title[Extension Theory]
{Operator Holes and extensions of sectorial operators
       and dual pairs of contractions}
\author{Mark~Malamud}
\address{Department of Mathematics \\
Donetsk National University \\
Universitetskaya str. 24 \\
83055 Donetsk \\
Ukraine}
\email{mdm@dc.donetsk.ua}

\keywords{Completion, shorted operator, generalized Schur complement,
selfadjoint contractive extension, nonnegative selfadjoint extension,
Friedrichs and Kre\u{\i}n-von Neumann extension}

\subjclass{Primary 47A57, 47B25; Secondary 47A55, 47B65.}

      \begin{abstract}
A description of the set of $m$-sectorial extensions of a dual pair
$\{A_1,A_2\}$ of nonnegative operators is obtained.
Some classes of nonaccretive extensions of the dual pair $\{A_1,A_2\}$
are described too. Both problems are reduced to similar problems for
a dual pair $\{T_1,T_2\}$ of nondensely defined symmetric contractions
$T_j=(I-A_j)(I+A_j)^{-1},\ j\in\{1,2\}$.
In turn these problems are reduced to the investigation of
the corresponding operator "holes".
A complete description of the set of all proper and improper extensions
of a nonnegative operator is obtained too.
    \end{abstract}





\maketitle


\bigskip

\section{Introduction} \label{intro}

In the theory of extensions of a nonnegative operator $A(\subset A^*)$ in
a Hilbert space $\gH$ to a selfadjoint or $m$-sectorial \cite{Ka}
operator
there are two
well-known approaches in which extensions $\tA\supset A$ in various
classes are described in diverse forms.
One of these, proposed by M.~G.~Krein in \cite{Kr1} (see also
\cite{AG, RSN}) uses the linear fractional transformation
$T_1=(I-A)(I+A)^{-1}$  to reduce the problem to the
description of various classes of extensions $T\supset T_1$
of a nondensely defined (on the subspace $\gH_1=(I+A)\gH$)
symmetric  contraction $T_1$.

The other approach to the description of proper extensions $\tA$ of an
operator $A>0$ was proposed by Vishik \cite{V} and Birman \cite{B}.
They associate with each extension $\tA\supset A$
(not necessarily selfadjoint) a "boundary" operator $B$
acting in an auxiliary space
$\cH\bigl(dim\cH=dim(A^*-i)\gH \bigr)$, and they describe the
properties of the extension $\tA=\tA_B$ in terms of the operator $B$,
i.e. essentially in terms of the boundary conditions if $A$
is a differential operator.
This approach was subsequently formalized in the concept of a
"boundary triplet' and was developed in later papers by many authors
(see for instance \cite{GG, DM1} and and references therein).

We remark that the methods used in these approaches
are essentially different,
as are the descriptions obtained with their help.

Recall that a closed densely defined operator $A$ in $\gH$ is
called sectorial with a half-angle $\ff\in(0,\pi/2]$ if
      \begin{equation}\label{1.1A}
 Re(Af, f)\ge  \ctg\ff \cdot |Im (Af,f)|,\qquad f\in\dom A.
     \end{equation}
It is called a maximal sectorial ($m$-sectorial) and is put in class
$S_{\gH}(\ff)$ if additionally $\rho(A)\not =\emptyset$.
If $\ff=\pi/2$ inequality \eqref{1.1A} turns into the inequality
$Re(A f,f)\ge 0$ and the class $S_{\gH}(\pi/2)$
is  the class of maximal accretive operators.  Denote also by
$S_{\gH}(0)$ the class of nonnegative selfadjoint operators
in $\gH$ and  note that  $S_{\gH}(0)=\cap_{\ff>0}S_{\gH}(\ff).$

In this paper we solve among others the following two problems.

{\bf Problem 1S.} Given a closed nonnegative symmetric operator
$A\ge 0$ in $\gH$. Describe the set $\Ext_A(\ff)$ of all proper and improper
$S_{\gH}(\ff)$-extensions of $A$ with $\ff\in[0,\pi/2]$.

{\bf Problem 2S.} Given a dual pair $\{A_1,A_2\}$ of closed
nonnegative symmetric operators in $\gH$.
Find necessary and sufficient conditions for
$\{A_1,A_2\}$ to admit an extension
${\wt A}\ (A_1\subset{\wt A}\subset A^*_2)$ of the class
$S_{\gH}(\ff)$ with $\ff\in[0,\pi/2]$ and describe the set
$\Ext_{\{A_1,A_2\}}(\ff)$ of  such extensions.

Note, that Problem 1S is solvable  for any $\ff\in [0, \pi/2].$
Indeed, it is known (see \cite{F, AG, RSN}) that
any symmetric
operator $A\ge 0$ admits a selfadjoint extension ${\wt A}\ge 0,$
say the Friedrichs extension $A_F.$
In other words, $\Ext_A(0)\not =\emptyset,$
hence $\Ext_A(\ff)\not =\emptyset$ for any
$\ff\in (0, \pi/2].$

A complete description of the set $\Ext_A(0)$
in terms of "boundary" operators
have been obtained in \cite{B} in the case of
a positive definite operator $A.$
The set $\Extp_A(\ff):=\Ext_{\{A,A\}}(\ff),\ $
$\ff\in [0,\pi/2],$ of all proper $m$-sectorial
extensions of an operator $A\ge 0$ with zero lower
bound was described via boundary  triplets and
Weyl functions in \cite{KM1} and \cite{DMT}.
Another description in the framework
of Krein's approach has been obtained in
\cite{AT1, AT2}.

On the other hand, even a solvability criterion
of Problem 2S was unknown.
We will show below that Problem 2S
is not necessary  solvable for any $\ff\in[0,\pi/2).$
It may even happen that it is solvable only with $\ff=\pi/2$.

We will also discuss  the  following more general problems.

{\bf Problem 3S.} Given a sectorial operator $A$ with
a half-angle $\ff_0\in [0,\pi/2).$
Describe the set $\Ext_A(\ff)$ of all
$S_{\gH}(\ff)$-extensions of $A$  with $\ff\ge\ff_0.$

{\bf  Problem 4S.}
Given a dual pair $\{A_1, A_2\}$ of sectorial operators in
$\gH$. Find necessary and sufficient conditions for
$\{A_1,A_2\}$ to have an extension
${\wt A}\ (A_1\subset{\wt A}\subset A^*_2)$
belonging to the class $S_{\gH}(\ff)$ with $\ff\ge\ff_0$ and
describe the set $\Ext_{\{A_1,A_2\}}(\ff)$ of all such extensions.

By the Kato-Schechter theorem (see \cite{Ka})
any sectorial operator $A$ obeying \eqref{1.1A}
with $\ff_0\in (0,\pi/2)$ admits
$m$-sectorial extension, say the Friedrichs
extension $A_F.$  In other words, $\Ext_A(\ff_0)\not =\emptyset,$
hence $\Ext_A(\ff)\not =\emptyset$ for $\ff\ge\ff_0,$
Thus,  Problem 3S is solvable for any $\ff\ge\ff_0.$

Note, that even a criterion of solvability of Problem 4S is unknown.

In accordance with  Krein's approach we consider a linear
fractional transformation
$T_1=(I-A)(I+A)^{-1}$ of a sectorial operator $A.$
It is clear that  $T_1$ is a nondensely defined contraction,
$T_1(\in[\gH_1,\gH])$, obeying the following condition
      \begin{equation}\label{1.2A}
\|T_1\sin\ff \pm i\cos\ff\cdot I\| \le 1 \qquad \dom T_1=\gH_1:=\ran(I+A).
    \end{equation}
We put an operator $T$ in the class $C_{\gH}(\ff)$ with
$\ff\in (0, \pi/2]$ if $\dom T=\gH$ and
inequality \eqref{1.2A} holds with $T$ in place of $T_1$.
Note that $C_{\gH}(\pi/2)$ is the class of all contractions in
$\gH$  and denote by $C_{\gH}(0)$ the class of all selfadjoint
contractions in $\gH$.

Now we can reformulate Problems 1S---4S in the following way.

{\bf Problem 1C.}
Given a nondendely defined symmetric operator $T_1\in[\gH_1,\gH].$
Describe the set of all  proper and improper
$C_{\gH}(\ff)$-extensions of $T$ with $\ff\ge\ff_0.$

{\bf Problem 2C.}
Given a dual pair $\{T_1,T_2\}$ of nondensely defined
symmetric contractions.
Find necessary and sufficient conditions for $\{T_1,T_2\}$
to admit an extension $T\in C_{\gH}(\ff)$ with
$\ff\ge\ff_0$ and  describe  the set
$\Ext_{\{T_1,T_2\}}(\ff)$ of all such extensions.

{\bf Problem 3C.}
Given a nondendely defined operator $T_1\in[\gH_1,\gH]$
obeying \eqref{1.2A}  with $\varphi=\varphi_0$.
Describe the set $\Ext_{T_1}(\ff)$ of all
$C_{\gH}(\varphi)$-extensions of $T$ with $\ff\ge\ff_0.$

{\bf Problem 4C.}
Given a dual pair $\{T_1,T_2\}$ of nondensely defined
contractions, obeying condition \eqref{1.2A} with
$\ff=\ff_0(\in[0,\pi/2])$.
Find necessary and sufficient conditions for $\{T_1,T_2\}$
to have an extension $T\in C_{\gH}(\ff)$
with $\ff\ge\ff_0$ and  describe
the set $\Ext_{\{T_1,T_2\}}(\ff)$ of all such extensions.






It is convenient to regard
Problems 1C and 3C as a problem on the "completion"
of a contractive operator matrix $T_1=\binom{T_{11}}{T_{21}}$ to form
a matrix $T=(T_{jk})^2_{j,k=1}$ which is connected in
a natural way with the problem of extending of a dual
pair of contractions  to operators in
various classes. A description is given in terms of
operator balls, "holes", and objects close to them.

Starting point of our investigation is a description of
the set of all contractive extensions of a
dual pair of contractions
$\bigl\{T_1=\binom{T_{11}}{T_{21}}, \ T_2=\binom{T^*_{11}}{T^*_{12}}\bigr\}$
or what is the same a description of all
"completions" of a matrix
     \begin {equation} \label{1.1}
  T_0
  = \begin{pmatrix} T_{11}& T_{12}\\ T_{21}&* \end{pmatrix}
  = \begin{pmatrix} T_{11}& D_{T_{11}^*}U\\ VD_{T_{11}}&* \end{pmatrix}
          \end{equation}
to form a contractive matrix $T=(T_{ij})^2_{i,j=1}$.

It has been shown in \cite{ArG, Da, DKW, ShY}
that all missing blocks $T_{22}$
in \eqref{1.1} form an operator
ball $B(-VT_{11}^*U; D_{V^*}, D_U):$
    \begin{equation} \label{1.2}
T_{22}=-VT^*_{11}U+D_{V^*}KD_U,\qquad  ||K||\le 1.
     \end{equation}
Our approach to  Problems 1C--4C is essentially based
on the solution to the following

    { \bf Problem 5.}
Given two operator balls
    \begin{equation*}
B(C_{\pm}; R^{\pm}_l, R^{\pm}_r)=
\{Z\in [\cH]:\  Z=C_{\pm} + R^{\pm}_l K R^{\pm}_r, \quad ||K||\le 1\}.
     \end{equation*}
Find a criterion for an operator "hole" ("loone")
     \begin{equation} \label{1.3}
L:=B(C_+; R^+_l, R^+_r) \cap B(C_-; R^-_l, R^+_r)
      \end{equation}
to be nonempty and  obtain a parametrization of $L.$

Problem 5 naturally arrises in diferent areas and is of interest
itself.
We will show here  that all Problems 1C-4C  are reduced to Problem 5.
Analysis of operator holes \eqref{1.3}
corresponding to Problems 1C-4C  shows that degree
of difficulty of any Problem jC with $j\in \{1,...,4\},$
can be characterized by means of the corresponding
radii $R_l^{\pm}$  and   $R_r^{\pm}.$
From this point of view Problem 2C with
$T_1=T_2 (\equ T_{11}=T^*_{11},\ T_{21}=T^*_{12})$
is the simplest one.  It is reduced to  Problem 5
with four equal radii $R_l^{\pm}=R_r^{\pm}=D_U$.
This problem is always solvable and it is equivalent
to a description  of the set $\Extp_{T_1}(\ff)$ of
proper $C(\ff)$-extensions of a symmetric  contraction $T_1,$
which has been  solved in \cite{AT1, AT2} by different  method.

Next, a solution to Problem 2C  is
equivalent to a description of missing blocks $T_{22}$
in  matrix \eqref{1.1} (with $T_{11}=T_{11}^*$)
such that $T=(T_{ij})\in C_{\gH}(\ff),$ that is
$T\sin\ff \pm i\cos\ff\cdot I\in C_{\gH}(\pi/2).$
Due to \eqref{1.2} this problem is
reduced to  Problem 5 with $R_l^+=R_l^-=D_{V^*}$ and
$R_r^+=R_r^-=D_U$. It is not always solvable in general
(see below).

Further, Problem 3C with $\ff_0>0$ is reduced
to Problem 5 (see \cite{M4} and Remark \ref{remProblem3})
with different left radii $R^+_l\not =R^-_l$ and equal right
radii $R^+_r=R^-_r,$ while  it is always solvable.

Finally, the most difficult Problem 4C with $\ff_0>0$ is reduced to
Problem 5 with different left radii $R^+_l\not =R^-_l$ and
different right radii $R^+_r\not =R^-_r$
(see Proposition \ref{prproblem4}).

I don't know a criterion of solvability
of Problem 5 if either
$R^+_l\not =R^-_l$ or $R^+_r\not =R^-_r,$
while a parametrization of the hole $L$ can be
easily obtained if at least one of its elements
is known (see \cite{KM1, M4}).
However a  solution to Problem 5 with
$R^+_l=R^-_l$  and $R^+_r=R^-_r$ is rather simple
and is contained in Lemma \ref{lem4.2}.


The paper is organized as follows.

In Section 2 we summarize some definitions and statements which are
necessary in the sequel.

In Section 3 we present a solution to Problem 2C
(see Theorem \ref{th4.1})  based on
Lemma \ref{lem4.2} on a parametrization of an operator
hole \eqref{1.3} with $R^+_l=R^-_l$  and $R^+_r=R^-_r.$
It is worth to note that though
$\Ext_{\{T_1,T_2\}}(\pi/2)\not =\emptyset,$
it may happen that
$\Ext_{\{T_1,T_2\}}(\ff)=\emptyset$ for any
$\ff\in[\ff_0, \pi/2)$.
The solvability of Problem 2C depends on the operator
$$
  Q_0 =  D_V^{-1}(I-VU)D_U^{-1}.
$$
More precisely, $\Ext_{\{T_1,T_2\}}(\ff)\not =\emptyset$
if and only if $\ff\in [\ff_1, \pi/2]$
where $\ff_1=\arccos(\|Q_0\|^{-1})$.
In particular,  $\Ext_{\{T_1,T_2\}}(\ff)=\emptyset$
for any $\ff\in (0, \pi/2)$ if and only if $Q_0$ is unbounded.

Further, in Section 3 we present a description of the set
${\Ext}_{T_1}(\ff)$ of all (proper and  improper)
extensions  of a symmetric contraction
$T_1\left(\in[\gH_1,\gH]\right)$ (see Theorem \ref{th4.2}).
This result gives a complete  solution to Problem 1C.

We also present here (see Propositions \ref{prop4.1A}
and \ref{prop4extreme}) a partial description of the
set $\Ext^e_{\{T_1,T_2\}}(\ff)$ of extreme points of the set
$\Ext_{\{T_1,T_2\}}(\ff)$.
It is interesting to note that even in a finite
dimensional case $(\dim\gH=n<{\infty})$
the set $C^e_{\gH}(\ff)$  of extreme points of the operator
loone $C_{\gH}(\ff)$ with $\ff\in(0,\pi/2)$
essentially differs from  the set $C^e_{\gH}(\pi/2)$ of extreme
points of the operator ball in ${\Bbb C}^n$.
Namely, though the set $C^e_{\gH}(\pi/2)$ consists of
unitary matrices, the set $C^e_{\gH}(\ff),\  \ff\in(0,\pi/2),$
in addition to normal matrices with "boundary spectrum"
contains continuum nonnormal matrices  with "nonboundary" spectrum.

Finally, in  Proposition \ref{prproblem4},
we discuss a reduction of Problem 4C  to Problem 5.

In Section 4  we investigate  noncontractive extensions
of  a dual pair $\{T_1,T_2\}$ of symmetric contractions.
Namely, we consider a (not necessary contractive)
extension $T_K$ of the form \eqref{1.1}, \eqref{1.2} and
calculate the Schur complement of any of the operators
      \begin{equation*}
G^\pm := (G_{ij}^\pm)_{i,j=1}^2 := I-T_KT_K^*\pm i\ctg\ff(T_K-T_K^*),
\qquad \ff\in (0, \pi/2].
      \end{equation*}
More precisely, assuming (for simplicity)
that $0\in\rho(G_{11})$ we prove
(see Theorem \ref{th4.3}) the following identities
     \begin{equation}\label{1.shur}
\sin^2\ff\cdot [G_{22}^\pm - G_{21}^\pm G_{11}^{-1}G_{12}^{\pm}]
  = D_U\cdot [I-(K^*\sin\ff\mp iQ^*)(K\sin\ff\pm iQ)]\cdot D_U,
         \end{equation}
where $Q$ is the closure of $Q_0$.

Using  \eqref{1.shur}  we describe the classes
$C_{\gH}(\ff;\kk^{\pm})$ of operators $T_K$ obeying conditions

\noindent
$\dim\ran(G^{\pm})_-=\kk^{\pm}$, where $\kk^{\pm}\in{\Bbb Z}_+$ and
$G_-$ stands for the "negative" part of the operator $G=G^*$.
Some applications of this result to the boundary value
problems can be found in \cite{M3}.
Moreover, formula \eqref{1.shur}
makes it possible to give another solution to Problem 2C
as well as to  obtain some complements to Theorem \ref{th4.1}.

In Section 5 we investigate completions of an incomplete matrix
$T'_0=
\begin{pmatrix}
T_{11}&*\\
0&T_{22}
\end{pmatrix}.$
Namely, in  Proposition \ref{pr4.2} we describe the set
of  some classes of noncontractive completions
of $T'_0$. This result complements and
generalizes the result of Nagy and Foias  \cite{SNF1}.

Moreover, in Proposition \ref{pr4.3}
we describe the sets of $C_{\gH}(\ff)$-completions of
$T'_0,$ giving an answer  to Yu. L. Shmul'yan's question.
This description is given in terms of operator holes.

Some results of the paper have been announced in \cite{M2}
and partially published (with proofs) in \cite{KM1}.

{\bf Notations.}  By  $\gH$ and $\cH$ we denote separable Hilbert
spaces; $[\gH_1,\gH_2]$ stands for the set of all bounded
linear operators from $\gH_1$ to $\gH_2;\  [\gH]:=[\gH,\gH];\  \cC(\gH)$
stands for the set of closed operators in $\gH$.
We denote by $\rho(T), \sigma(T)$ and $\sigma_{pp}(T)$ the resolvent set,
the spectrum and the purely point spectrum of
$T(\in\cC(\gH))$ respectively;
$\sigma_p(T)$ stands for the set of eigenvalues of $T;$
$\dom T$ and $\ran T$ stand for the domain of definition
and the range of the operator $T$ respectively.
As usual $E_T(\cdot)$ stands for the spectral measure
(resolution of the identity) of a self-adjoint
operator $T\in \cC(\gH);$ \  $T_-:=TE_T(0,\infty).$

\section{Preliminaries} \label{sect1}

\noindent{\bf 2.1. Dual pairs of contractions.}
We recall a definition of a dual pair of bounded operators.
    \begin{definition}\label{def2.1}
Let $\gH=\gH_1\oplus\gH_2=\gH_1'\oplus\gH_2'$ be
orthogonal decompositions of the Hilbert space $H$.
Operators $T_1\in[\gH_1,\gH],\
T_2\in[\gH_1',\gH]$ are said to form a dual pair
of bounded operators  if
    \begin{equation}\label{2.1}
(T_1f,g)=(f,T_2g),\qquad  f\in\gH_1,\    g\in \gH_1'.
     \end{equation}
An operator $T(\in[\gH])$ is termed an extension of the dual pair $\dup$ if
$$
T\lceil\gH_1=T_1 \quad \text{ and } \quad  T^*\lceil\gH_2=T_2.
$$
The set of all extensions of a dual pair $\dup$ is denoted
by $\Ext_{\{T_1, T_2\}}.$
   \end{definition}
When rewritten in the block-matrix representation with
respect to the pointed  out decompositions of
the space $\gH$, the operators $T_1$ and $T_2$ form a
dual pair if and only if
    \begin{equation}\label{2.2}
T_1=\binom{T_{11}}{T_{21}},\qquad T_2=\binom{T_{11}^*}{T_{21}'}
     \end{equation}
with $T_{11}\in[\gH_1,\gH_1'],\ T_{21}\in[\gH_1,\gH_2'],\
T_{21}'\in[\gH_1', \gH_2]$.

Setting $T_{12}=(T_{21}')^*$, an extension $T$ of the DP $\dup$ can be
rewritten in the form
    \begin{equation}\label{2.3}
T=\begin{pmatrix}
T_{11}&T_{12}\\
T_{21}&T_{22}
\end{pmatrix} \quad \text{with }\quad T_{22}\in [\gH_2,\gH_2'] .
     \end{equation}
In this case the problem of description of a certain class $X$
of extensions of the dual pair $\dup$ is equivalent
to the problem of completing an  incomplete block-matrix
$\begin{pmatrix}
T_{11}&T_{12}\\T_{21}&*
\end{pmatrix}$
with  respect to the matrix $T$ of the form \eqref{2.3}
and such that $T\in X$.

In  what follows  we consider contractive extensions of a
dual pair of contractions  $\dup$.
The union of all such extensions will be denoted by
$\Ext_{\{T_1,T_2\}}(\pi/2)$.

The set $\Ext_{\{T_1,T_2\}}(\pi/2)$ turns out to be an operator
ball in the sence of the following definition.
     \begin{definition}\label{def2.2}
The totality of the operators  $Z\in [\gH]$ of the form
    \begin{equation}\label{2.4}
Z=C_0+R_lKR_r,\qquad \|K\|\le 1
    \end{equation}
is referred to as an operator ball $B(C_0;R_l,R_r)$.

Here $C_0$ is called the center of the ball, and
 $R_l=R_l^*\ge 0$ and $R_r=R_r^*\ge 0$
are called left and right radii respectively.
       \end{definition}
We will use the following simple and known result.
    \begin{lemma}\label{lem2.2}\cite{GL}
Let $Q_j\in[\gH],\ j\in\{1,2,3\},\ Q_3=Q_3^*,\ Q_1>0$
and $0\in\rho(Q_1).$
Then the iequality
   \begin{equation} \label{2.6}
Z^*Q_1Z+Z^*Q_2+Q_2^*Z+Q_3\le 0
    \end{equation}
has a solution if and only if
    \begin{equation} \label{2.7}
Q_2^*Q_1^{-1}Q_2-Q_3\ge 0.
    \end{equation}
Under this condition the set of the solutions of the
inequality \eqref{2.6} makes up an operator
ball $B(C_0;R_l;R_r)$ of the form \eqref{2.4} with
   \begin{equation} \label{2.8}
C_0=-Q_1^{-1}Q_2,\quad R_l=Q_1^{-1/2}\quad \text{and}\quad
   R_r=(Q_2^*Q_1^{-1}Q_2-Q_3)^{1/2}.
    \end{equation}
\end{lemma}

   {\bf 2.2.} The operators $T_1$ and $T_2$ of the
form \eqref{2.2} are  contractive if and only if
     \begin{eqnarray}
&T_{11}^*T_{11}&+T_{21}^*T_{21}\le I\equ T_{21}^*T_{21}\le D_{T_{11}}:=
I-T_{11}^*T_{11},\\
\label{2.12}
&T_{11}T_{11}^*&+T_{12}T_{12}^*\le I\equ T_{12}T_{12}^*\le D_{T_{11}^*}:=
I-T_{11}T_{11}^*.
      \end{eqnarray}
It is known  (and it is obvious) that
these relations are equivalent to the  following ones
    \begin{equation}\label{2.13}
T_{21}=VD_{T_{11}},\qquad T_{12}=D_{T_{11}^*}U
     \end{equation}
with contractions $V$ and
$U\ (V\in[\gH_1,\gH_2'],\ U\in[\gH_2,\gH_1']),$
which  are uniquelly  determined
provided that $\ker V\supset \ker D_{T_{11}}$ and
$\ker U^*\supset \ker D_{T_{11}^*}$,
that is $V^*=D_{T_{11}}^{-1}T_{21}^*$
and  $U=D_{T_{11}^*}^{-1}T_{12}$.
A complete description of the set
 $\Ext_{\{T_1,T_2\}}(\pi/2)$  is contained
in the  following theorem.
We will  essentially use it in the sequel.
         \begin{theorem}\label{th2.1}(\cite{ArG,Da, DKW,ShY}).
Let $\dup$ be a dual
pair of contractions,
    \begin{equation}\label{2.14}
T_1=\binom{T_{11}}{T_{21}}=\binom{T_{11}}{VD_{T_{11}}},\quad
T_2=\binom{T_{11}^*}{T_{12}}=\binom{T_{11}^*}{U^*D_{T_{11}^*}},
     \end{equation}
${\cH}_1={\overline\ran}(D_U)$ and
${\cH}_2={\overline\ran}(D_{V^*})$. Then
the formula
      \begin{equation}\label{2.15}
T:=T_K=
\begin{pmatrix} T_{11}&D_{T^*_{11}}U\\
VD_{T_{11}}&T_{22}
\end{pmatrix},\qquad
T_{22}=-VT_{11}^*U+D_{V^*}KD_U,
       \end{equation}
establishes a bijective correspondence between all
contractive extensions
$T:=T_K=(T_{ij})\in\Ext_{\{T_1,T_2\}}(\pi/2)$ and
all contractions  $K\in[{\cH}_1,{\cH}_2]$.

Thus, the set $\Ext_{\{T_1,T_2\}}(\pi/2)$ forms an
operator ball  $B(C_0;R_l,R_r)$ with  the center
$C_0=-VT_{11}^*U$  and left and right  radii
$R_l=D_{V^*}$ and  $R_r=D_U$ respectively.
      \end{theorem}
      \begin{remark}\label{rem2.2}
Let us make some historical remarks concerning Theorem
\ref{th2.1}.
The case $T_{11}=T_{11}^*,\ T_{12}=T_{21}^*$
was considered by M.G. Krein \cite{Kr1} in connection
with selfadjoint extensions of positive
unbounded operators while his description of the class
$Ext_{\{T_1,T_1\}}(0)$ differs
from that followed from Theorem \ref{th2.1}.
The existence of
contractive extensions of a DPC $\dup$ (that is the fact
$\Ext_{\{T_1,T_2\}}(\pi/2)\not=\eset$) was first established by
B.S. Nagy and C. Foias \cite{SNF2}, p.190,
by means of a corresponding
generalization of the Krein's method  \cite{Kr1}.
Note also that the claim
$\Ext_{\{T_1,T_2\}}(\pi/2)\not=\eset$
is implicitly  contained in \cite{Phi1, Phi2}.
Another proof of the existence part
of Theorem \ref{th2.1} has also been obtained by
S. Parrot \cite{Par}.

The complete description of the set $\Ext_{\{T_1,T_2\}}(\pi/2)$, i.e.,
Theorem \ref{th2.1}, was obtained in
\cite{ArG, Da, DKW, ShY}.
In the special case $(T_{21}=0)$
Theorem \ref{th2.1} has been
obtained by B.Sz.-Nagy and C.Foias \cite{SNF1}
much earlier.
Several other proofs  of Theorem \ref{th2.1},
based on different ideas, can also  be found in \cite{FF, KM1, M4}.
In particular the proof of  C. Foias and A.E. Frazho \cite{FF}
is based on Redheffer's  products,
the author's proof in \cite{M4} is based on Lemma \ref{lem2.2}.
  \end{remark}

{\bf 2.3. Extreme points of the unit ball.}

Recall the following

     \begin{definition} \label{def2.3}
Let $G$ be a closed convex set in a Banach space $X$.
A point $f\in G$ is called an extreme point of $G$ if it does not admit
a represetation
$f=f_1+(1-t)f_2$ with $f_1,f_2\in G, f_1\not =f_2$ and
$t\in(0,1)$.
     \end{definition}
Let $\gH_1,\gH_2$ be Hilbert spaces. An operator
$T\bigl(\in[\gH_1,\gH_2]\bigr)$ is called a partial isometry if
$T^* T=P$, where $P$ is an orthoprojection in $\gH_1$. An operator
$T\bigl(\in[\gH_1, \gH_2]\bigr)$ is a maximal partial isometry if either
$T$ or $T^*$ is an isometry from $\gH_1$ to $\gH_2$, that is, if either
$T^* T=I_{\gH_1}$ or $T T^*=I_{\gH_2}$.

In the sequel we need the following result which is well
known in the case $\gH_1=\gH_2$ (see \cite{Hal}).
The general case can be easily derived from the known one.
       \begin{proposition}\label{prop2.1}
 The set  $\gR^e_1$ of extreme points of the unit ball
$\gR_1:=\{T:T\in[\gH_1,\gH_2], \  \|T\|\le 1\}$ in
$[\gH_1,\gH_2]$ consists
of maximal partial isometries from $\gH_1$ to $\gH_2$.
     \end{proposition}

{\bf 2.4. Sectorial operators and  $C(\ff)$-contractions.}

\begin{definition}\label{def2.4}\cite{Ka}.
A closed linear operator  $A$ in a Hilbert space $\gH$
is called sectorial with vertex zero and half-angle
$\ff\in (0;\pi/2)$ if its numerical range
is contained in sector
$G_{\ff}=\{z\in\Bbb C:\  |\arg z|\le \ff<\pi/2\},$
that is
    \begin{equation} \label{3.1}
\ctg\ff\cdot|Im(Af, f)|\le \text{Re}(Af, f),\qquad f\in \dom A.
    \end{equation}
If in addition $A$ has no sectorial extensions
$(\Llr \rho(A)\not=\emptyset)$
it is called a $m$-sectorial operator
and is put in the class $S_{\gH}(\ff).$

Further by $S_{\gH}(\pi/2)$ we denote the class of
$m$-accretive operators in $\gH,$ i.e.
$A\in S_{\gH}(\pi/2)$ if
Re$(Af, f)\ge 0$ for all $f\in \dom A$
and  $\rho(A)\not=\emptyset.$

Finally, $S_{\gH}(0)$ stands for the set of all
nonnegative selfadjoint linear operators  in $\gH.$
          \end{definition}

Following  \cite{KrS} an operator $A$
with $\overline{\dom (A)}=\gH$ is called regularly
dissipative if $-A$ is $m$-sectorial.

Let $A$ be a closed sectorial closed  operator in $\gH.$
In the framework of  the approach accepted in this
paper with each $A$  it is connected a linear transformation
   \begin{equation} \label{3.2}
T_1=X(A):=-I+2(I+A)^{-1},
   \end{equation}
being a contraction with a nondense in $\gH$
domain of the definition
$\gH_1:=\dom(T_1)=(I+A)\dom A.$
In so doing  condition \eqref{3.1}
is transformed to the following one
    \begin{equation} \label{3.3}
2\ctg\ff\cdot |Im(T_1f,f)|\le ((I-T_1^*T_1)f,f)=\|D_{T_1}f\|^2,
\qquad f\in\dom T_1.
     \end{equation}

The following definition naturally arises from what has been said.
      \begin{definition}\label{def3.2}
We put an operator $T\in [\gH]$ in the class
$C_{\gH}(\ff)$ if
   \begin{equation} \label{3.4}
\|T\sin\ff\pm i\cos\ff\cdot I\|\le 1,\qquad (\ff\in(0,\pi/2])
   \end{equation}
and in the class $C_{\gH}(0)$ if $T=T^*$ and $\|T\|\le 1.$
        \end{definition}
It is clear that if $T\in C_{\gH}(\ff)$, then $\sigma(T)\subset L_{\ff}$
where
    \begin{equation}\label{3.4hole}
L_{\ff}:=\{z\in{\Bbb D}:\  |z\sin\ff\pm i\cos\ff|\le 1\}.
      \end{equation}
    \begin{lemma}\label{lem3.1}
Let $T\in [\gH],\quad \ff\in(0,\pi/2]$,
Then the following
properties of the operator $T$ are equivalent:
      \begin{eqnarray*} 
&(1)&  \ T\in C_{\gH}(\ff);  \\
&(2)&  \ 2\ctg\ff|Im(Tf,f)|\le \|D_Tf\|^2, \qquad f\in\gH; \\
&(3)& \ 2\ctg\ff|(T_If,g)|\le \|D_Tf\|\cdot\|D_Tg\|, \qquad  f,g\in\gH.
    \end{eqnarray*}
     \end{lemma}
It follows from Lemma \ref{lem3.1} that
      \begin{equation}\label{3.5}
C_{\gH}(0)=\cap_{\ff\in(0,\pi/2)}C_{\gH}(\ff).
     \end{equation}
     \begin{definition}\label{def3.2A}
Let $A$ be a closed sectorial operator in  $\gH$
with  vertex zero and half-angle
$\ff_0\in [0;\pi/2)$ and let $T_1$ be a
contraction obeying  \eqref{3.3}  with $\ff=\ff_0.$
Denote by ${\Ext}_A(\ff)$  the class of all
$m$-sectorial  extensions  ${\wt A}(\in \cC(\gH))$
of $A$ with  vertex zero and the half-angle
$\ff\in [\ff_0;\pi/2]$
and by  ${\Ext}_{T_1}(\ff)$ the class of
all extensions  $T\in[\gH]$ of $T_1$ obeying  \eqref{3.3}
with $\ff\in [\ff_0;\pi/2].$
     \end{definition}

    \begin{lemma}\label{lem3.2}
Let $A$ be a densely defined sectorial operator in $\gH$ with
a semiangle $\ff$.
The linear fractional transformation \eqref{3.2}
establishes the bijective  correspondence
     \begin{equation}
{\wt A}\to T=X({\wt A})=-I+2(I+{\wt A}^{-1}, \quad
T\to {\wt A}=X^{-1}(T)=-I+2(I+T)^{-1},
     \end{equation}
between the set  ${\Ext}_A(\ff)$  and the subset
${\Ext}'_{T_1}(\ff)=\{T:\ T\in {\Ext}_{T_1}(\ff),\  -1\in\rho(T)\}$
of  ${\Ext}_{T_1}(\ff).$
         \end{lemma}
If $A$ is a nondensely defined sectorial operator,
then the set ${\Ext}_A(\ff)$
contains  $m$-sectorial linear  relations too.
Lemma \ref{lem3.2} remains valid in this case if we replace
${\Ext}'_{T_1}(\ff)$ by ${\Ext}_{T_1}(\ff)$.

\section{Some classes of contractive extensions of
dual pairs  of Hermitian  contractions} \label{sect2}

{\bf 3.1. A parametrization of the operator loone in  the
special case.}

In this subsection we present an elementary result
(see Lemma \ref{lem4.2}) on parametrization of an operator hole
$L=B_1\cap B_2$ in the case of operator balls
$B_1=B(Z_1;R_l,R_r)$ and $B_2=B(Z_2;R_l,R_r)$
with equal left radii and right radii.
This lemma  gives  a partial solution
to Problem 5 mentioned in the Introduction.

We start with the following simple lemma.
      \begin{lemma}\label{lem4.1}
Let $R_1=R_1^*\ge 0,\ R_2=R_2^*\ge 0,\ R_i\in [\gH], \ j\in\{1,2\},$
and let $A\in [\gH]$. Then the following conditions are equivalent:
    \begin{eqnarray} \label{4.1}
&(i) &  A=R_2BR_1, \qquad B\in C_{\gH}(\pi/2);  \\
\label{4.2}
&(ii) &\ 2|(Af,g)|\le (R_1^2f,f) + (R_2^2g,g),\qquad  f,g\in \gH.
     \end{eqnarray}
         \end{lemma}
      \begin{proof}
Inequality \eqref{4.2} is equivalent to the inequality
     \begin{equation} \label{4.3}
|(Af,g)|\le \|R_1f\|\cdot\|R_2g\|,\qquad  f,g\in \gH.
      \end{equation}
Hence $\ker A \supset \ker R_1 $ and  $\ker A^*\supset \ker R_2$.
Letting
$R_1f=:h_1,\ R_2g=:h_2$ one rewrites \eqref{4.3} in the form
    \begin{equation*}
|(AR_1^{-1}h_1, R_2^{-1}h_2)|\le \|h_1\|\cdot \|h_2\|,\qquad
h_1\in \ran R_1,\ h_2\in \ran R_2.
     \end{equation*}
It follows that  the bylinear form
$(AR_1^{-1}h_1, R_2^{-1}h_2)$
may be continually extended to a bounded bylinear form
$\gt(h_1,h_2)$ on $\cH_1\times \cH_2$, with
$\cH_j=\overline{\ran}R_j,\ j\in \{1,2\}.$
Hence
$\gt(h_1, h_2)=B(h_1, h_2),$
where $B\in [\cH_1,\cH_2]$ and $\|B\|\le 1$. Since
$(AR_1^{-1}h_1, R_2^{-1}h_2)=(Bh_1, h_2),\  h_j\in\cH_j,$
\ $j\in\{1,2\},$  then  $AR^{-1}h_1 \in \dom(R_2^{-1})$ and
$R_2^{-1}AR_1^{-1}h_1 = Bh_1,\quad h_1\in \dom(R_1^{-1})$.
Thus, $B$ is the closure of the operator
$R_2^{-1}AR_1^{-1}.$
Hence  $A=R_2BR_1$ and the implication
(ii)$\then$(i) is proved.

The converse implication (i) $\then$(ii) is clear.
          \end{proof}

The following statement easily follows from
Lemma \ref{lem4.1}.
       \begin{lemma}\label{lem4.1A}
Let $R_j=R^*_j\ge 0, \  R_j\in[\gH],\  j\in\{1,2\}$
and $A\in[\gH]$. Suppose additionally that
$\ker R_1 = \ker R_2$ and
$\cH:=\gH\ominus \ker R_1$.
Then $A$ admits a representation $A=R_2 K R_1$ with
$K\in C_{\cH}(\varphi), \  \varphi\in[0,\pi/2]$,
if and only if
      \begin{equation}\label{2.4A}
2Re\bigl((\sin\varphi\cdot A\pm i \cos\varphi\cdot R_2\cdot R_1))f,g\bigr)
\le(R_1^2 f,f)+(R^2_2 g,g),\qquad  f,g\in\cH.
     \end{equation}
         \end{lemma}
      \begin{proof}
Necessity is immediately implied by Lemma \ref{lem4.1}.

Sufficiency. Suppose that \eqref{2.4A} is satisfied. Assume that
$\varphi>0$, since the case $\varphi=0$ is trivial.
Then \eqref{2.4A} yields
       \begin{equation*}
2|(A f,g)|\le(sin\varphi)^{-1}\cdot(\|R_1 f\|^2+\|R_2 g\|^2), \qquad
f,g\subset\gH.
    \end{equation*}
By Lemma \ref{lem4.1} $A$ admits a representation $A=R_2 K R_1$ with
$K\in [\cH], \|K\|\le 1/sin\varphi$.
Substituting this expression for
$A$ in \eqref{2.4A} we get the required.
      \end{proof}

The following lemma, being a partial solution
to  Problem 5, gives a parametrization of the
operator loone $L:=B_1\cap B_2$
in the case of operator balls  $B_1$ and $B_2$ in $[\gH]$
with equal left and right radii and,
in particular, it gives a criterion of nonemptyness
of the loone $L.$

           \begin{lemma}\label{lem4.2}
Let $B_1=B(C_1;R_l,R_r)$ and $B_2=B_2(C_2;R_l,R_r)$
be two operator balls in $[\gH]$ with equal left and right radii and
$\cH_1:=\gH\ominus \ker R_r,\   \cH_2:=\gH\ominus \ker R_l.$
Then

(i) their  intersection $L:=B_1\cap B_2$  is nonempty if and only
if one of the  following  (equivalent)  conditions
is satisfied:
      \begin{eqnarray} \label{4.4}
&(a)& \ |(C_1-C_2)f,g)|\le 2^{-1}[(R_l^2f,f)+(R_r^2g,g)];  \\
\label{4.5}
&(b)&  \ |(C_1-C_2)f,g)|\le \|R_lf\|\cdot \|R_rg\|,\qquad  f,g\in \gH; \\
\label{4.6}
&(c)&  \ C_1-C_2 = 2R_lQR_r  \quad  \text{with }\quad
Q\in C(\pi/2),\   Q\in [\cH_1,\cH_2];
      \end{eqnarray}
      \begin{multline*}
(c') \quad \text{the operator } \quad
Q_0:=2^{-1}R_l^{-1}(C_1-C_2)R_r^{-1}
\quad \text{ is bounded and its closure } \\
Q:=\bar Q_0 (\in [\cH_1,\cH_2])\
\text{ is a contraction}.
       \end{multline*}

(ii)\
If any of the conditions (a), (b), (c) is satisfied, then the operator
loone $L$  admits the following  parameter representation
    \begin{equation} \label{4.7}
T\in L=B_1\cap B_2\equ T=T_K:=2^{-1}(C_1+C_2)+R_lKR_r  \quad \text{with}\quad
K\pm Q\in C(\pi/2).
       \end{equation}

(iii)\  $L$ consists of one element, $L=\{2^{-1}(C_1+C_2)\},$
if and only if at least one of the following three conditions holds
$$
(1) \quad  R_l=0; \quad  (2)\quad   R_r=0; \quad
(3)\    Q \ \text{is a maximal partial isometry}.
$$
        \end{lemma}
      \begin{proof}
(i),\ (ii).   Equivalence of the conditions (a)-(c) is implied by
Lemma \ref{lem4.1}.
It remains to show, for example, that equality
\eqref{4.6} is equaivalent to the condition
$L\not=\eset$.  Let $T\in L=B_1\cap B_2$,  that is
     \begin{equation*}
T=C_1+R_lK_1R_r=C_2+R_lK_2R_r.
     \end{equation*}
Then setting
$$
K:=2^{-1}(K_1+K_2)(\in C(\pi/2))\quad \text{and} \quad
 Q:=2^{-1}(K_2-K_1)(\in C(\pi/2)),
$$
we deduce
     \begin{equation*}
T=2^{-1}(C_1+C_2)+R_lKR_r, \quad \text{where}\quad   K\pm Q\in C(\pi/2).
       \end{equation*}
Thus, conditions \eqref{4.6} and \eqref{4.7} are satisfied.

Conversely, suppose that  \eqref{4.6} is valid.
Then setting  $K_1:=-Q,\ K_2 :=+Q$  and  $T:=C_1+R_lK_1R_r$
we get
   \begin{equation*}
T=C_1+R_lK_1R_r = C_2+R_lK_2R_r.
    \end{equation*}
Hence    $T \in L= B_1\cap B_2.$

(iii)
Let $T\in L$ and $T\not =2^{-1}(C_1+C_2)$. Then according
to \eqref{4.7}
$R_l\not =0,\ R_r\not =0$ and there exists
$K\in C(\pi/2)\setminus\{0\}$ such that
$K_{\pm}:=Q\pm K\in C(\pi/2)$. Hence
$Q=(K_++K_-)/2\not =K_+\not =K_-$
is not an extreme point of the unit ball in
$[\frak H_1,\frak H_2]$. By Proposition \ref{prop2.1}
$Q$ is not a maximal isometry.

Conversely, suppose that $R_l\not =0,\ R_r\not =0$ and $Q$
is not a maximal  isometry. By Proposition \ref{prop2.1}
$Q=(K_++K_-)/2$ where
$K_{\pm}\in[\cH_1, \cH_2],\  K_{\pm}\in C(\pi/2)$ and
$K_+\not =K_-$.
Setting $K_1:=(K_+-K_-)/2$ and $K_2:=-K_1$ we easily get that
$K_j\pm Q\in C(\pi/2), \ j\in \{1,2\}$.
Hence by \eqref{4.7}
$T_{K_j}:=2^{-1}(C_1+C_2) + R_l K_j R_r\in L,\ j\in\{1,2\}.$
Since  $T_{K_j}\not =2^{-1}(C_1+C_2),$
we get the required.
           \end{proof}

{\bf 3.2. A description of the class of  $C_\gH(\ff)$-extensions
of a dual pair of symmetric contractions.}
In this subsection we present  a solution to the Problem 2C with
$\varphi_0=0$.

More precisely, let $\dup$ be a dual pair of symmetric contractions in
$\gH=\gH_1\oplus\gH_2=\gH_1'\oplus\gH_2'.$
Due to \eqref{2.14}  the operators $T_1$ and $T_2$
admit the following block-matrix representations
     \begin{equation}\label{4.8}
T_1=\binom{T_{11}}{T_{21}}=\binom{T_{11}}{VD_{T_{11}}},\qquad
T_2=\binom{T_{11}}{T_{12}^*}=\binom {T_{11}}{U^*D_{T_{11}}},
    \end{equation}
since  $T_{11}=T_{11}^*$.  In particular, in this case
$\gH_1=\gH_1'$ and $\gH_2=\gH_2'.$

Let
    \begin{equation}\label{4.8A}
\Ext_{\{T_1,T_2\}}(\varphi):=\Ext_{\{T_1,T_2\}}\cap C_\gH(\ff),\qquad
\ff\in [0,\pi/2],
     \end{equation}
stand for the set of $C_\gH(\ff)$-extensions of
the dual pair of symmetric contractions $\{T_1,T_2\}.$

By  Theorem \ref{th2.1}  $\Ext_{\{T_1,T_2\}}(\pi/2)\not =\emptyset,$
that is there  always exists an
extension  $T\in C_\gH(\pi/2)$ of the dual pair
$\dup.$  
It turns out that it is not  the case for the classes
$C_\gH(\ff)$ and $C_{\gH}(\ff;\kk)$ with $\ff<\pi/2$.
The solvability of both  problems  depends  on the
properties of the operator
    \begin{equation} \label{4.9}
Q_0:= D_{V^*}^{-1}(I-VU)D_U^{-1}.
     \end{equation}
Moreover, we show that if the operator $Q_0$ is unbounded
the Problem 2C  has a solution only with
$\varphi_1=\pi/2$, that is,
$\Ext_{\{T_1,T_2\}}(\varphi)\not =\emptyset$ iff
$\ff=\ff_1=\pi/2$.

          \begin{theorem}\label{th4.1}
Let $\dup$ be a dual pair of symmetric contractions
in $\gH=\gH_1\oplus \gH_2$  of the form \eqref{4.8},
$\ff\in (0,\pi/2]$
and let $\cH_1:=\overline{\ran}(D_U),\ \cH_2:=
\overline{\ran}(D_{V^*})$.
Then:

  (i)
the set  $\Ext_{\{T_1,T_2\}}(\ff)$   is nonempty if and only if
$\ff\in[\ff_1,\pi/2]$, where
      \begin{equation}\label{4.10}
\ff_1 := \arccos\bigl(\|D^{-1}_{V^*}(I-V U)D^{-1}_U\|^{-1}\bigr);
      \end{equation}

\noindent

  (ii)
for any  $\ff\in[\ff_1,\pi/2]$
the following equivalence holds: 
      \begin{multline} \label{4.11}
    T_K
    = \begin{pmatrix}
        T_{11} & T_{12} \\
        T_{21} & T_{22}
      \end{pmatrix}
    \in  \Ext_{\{T_1,T_2\}}(\ff)     \Llr  \\
    T_{22} = -VT_{11}U+D_{V^*}KD_U, \quad K\cdot\sin\ff\pm
iQ\cos\ff\in C(\pi/2),
          \end{multline}
where  $Q := {\overline Q}_0(\in [\cH_1,\cH_2])$ is the closure of
the operator $Q_0$ of the form \eqref{4.9}
and $K\in [\cH_1,\cH_2]$;

(iii) the set $\Ext_{\{T_1,T_2\}}(\ff)$
consists of one element
if at least one of the following three conditions is satisfied:

(a) $D_{V^*}=0;$  \qquad (b)\  $D_U=0;\qquad$
(c)\   $Q\cos\ff$  is a maximal partial isometry.
          \end{theorem}
       \begin{proof}
(i)  Let $T=(T_{ij})_{i,j=1}^2$ be a contractive extension
of the dual pair $\{T_1, T_2\}.$
Suppose for the begining that $\ff>0.$
In this case  the inclusion $T\in C_\gH(\ff)$  means that
$$
T_{\pm}:=T\sin\ff\pm i\cos\ff\cdot I\in C_\gN(\pi/2),
$$
that is,
  \begin{equation*}            
    T_{\pm} :=
    \begin{pmatrix}
      B_{11}^{\pm} & B_{12} \\
      B_{21}       & B_{2}^{\pm}
    \end{pmatrix}
    :=
    \begin{pmatrix}
      \sin\ff\cdot T_{11}\pm i\cos\ff\cdot I & T_{21}\sin\ff \\
      T_{21}\sin\ff & \sin\ff\cdot T_{22}\pm i\cos\ff\cdot I
    \end{pmatrix}
    \in C_\gH(\pi/2).
  \end{equation*}          
It is easily  seen that
  \begin{equation} \label{4.12}
    \begin{split} 
    D_{B_{11}^+}^2=I-(\sin\ff\cdot T_{11}-i\cos \ff\cdot I)(\sin\ff
\cdot T_{11}+i\cos \ff\cdot I)=\sin^2\ff\cdot D_{T_{11}}^2,  \\
    D_{B_{11}^-}^2=I-(\sin\ff\cdot T_{11}+i\cos \ff\cdot I)
(\sin\ff \cdot T_{11}-i\cos\ff\cdot I)=\sin^2\ff\cdot D_{T_{11}}^2.
    \end{split}
  \end{equation}
  Therefore $D_{B_{11}^+}=D_{B_{11}^-}=\sin\ff\cdot D_{T_{11}}$
  and consequently
\begin{equation*}
  \begin{split}
    B_{12}
      &= \sin\ff T_{12}
      = \sin\ff\cdot D_{T_{11}}U
      = D_{B_{11}^{\pm}}U, \\
    B_{21}
      &= \sin\ff T_{21}
      = \sin\ff\cdot VD_{T_{11}}
      = VD_{B_{11}^{\pm}}.
  \end{split}
      \end{equation*}
  Thus the contractions $T_{\pm}$ have the form
  \begin{equation} \label{4.13}
    T_{\pm} =
    \begin{pmatrix}
      B_{11}^{\pm} & B_{12} \\
      B_{21}       & B_{22}^{\pm}
    \end{pmatrix}
    =
    \begin{pmatrix}
      B_{11}^{\pm}      & D_{B_{11}^{\pm}}U \\
      VD_{B_{11}^{\pm}} & B_{22}^{\pm}
    \end{pmatrix}
  \end{equation}
with $B_{11}^{\pm}=T_{11}\sin\ff\pm i\cos\ff\cdot I$ and
$B_{22}^{\pm}=T_{22}\sin\ff\pm i\cos\ff\cdot I$.
According to Theorem \ref{th2.1} the equivalences
     \begin{equation} \label{4.14}
T_{\pm}\in C(\pi/2)\equ B_{22}^+=T_{22}\sin\ff\pm i\cos \ff\cdot I=C_{\pm}
+D_{V^*}K_{\pm}D_U
      \end{equation}
hold true with some contractions $K_{\pm}$ and operators $C_{\pm}$
defined by
     \begin{equation} \label{4.15}
C_{\pm}:=-V(B_{11}^{\pm})^*U=-V(\sin\ff\cdot T_{11}^*\mp i\cos \ff)U.
      \end{equation}
Setting
     \begin{equation} \label{4.16}
C_1:=C_+-i\cos\ff\cdot I,\qquad C_2:=C_-+i\cos \ff\cdot I,
      \end{equation}
we rewrite equivalences \eqref{4.14} in the form
      \begin{equation} \label{4.17}
T\in C_{\gH}(\ff)\equ T_{22}\sin\ff =
C_1+D_{V^*}K_{+}D_U=C_2+D_{V^*}K\_D_U, \quad
K_{\pm}\in C(\pi/2).
       \end{equation}
Thus, $T\in C_\gH(\ff)$ if and only if the operator
$\sin\ff\cdot T_{22}$ belongs to the intersection
of the operator balls $B_1:=B(C_1;D_{V^*};D_U)$ and
$B_2:=B(C_2;D_{V^*},D_U)$, that is,
$\sin\ff\cdot T_{22}\in L:=B_1\cap B_2$.
By Lemma \ref{lem4.2} with account of
\eqref{4.15} and \eqref{4.16} the
condition $L= B_1\cap B_2\not=\eset$ amounts to saying that the operator
       \begin{equation} \label{4.18}
C_0:=2^{-1}(C_1-C_2)=-i\cos\ff(I-VU)
        \end{equation}
admits the representation $C_0=D_{V^*}KD_U$ with $K\in C(\pi/2)$
or, what is the same,  the operator
$Q_0\cdot \cos\ff$  is contractive where the operator
$Q_0$ is  of the form \eqref{4.9}.
This proves the first assertion.

(ii)  Suppose that condition  \eqref{4.10}
is satisfied.
To obtain  a parametrization of the hole $L$
we note that by \eqref{4.15} and  \eqref{4.16}
$$
2^{-1}(C_1+C_2) =-\sin\ff\cdot V T^*_{11}U.
$$
Now Lemma \ref{lem4.2}  yields the equivalence
      \begin{equation}\label{4.18a}
T\in L (=B_1\cap B_2)  \Llr  \sin\ff\cdot T_{22}=
-\sin\ff\cdot VT^*_{11}U + D_{V^*}{\widetilde K}D_U,
    \end{equation}
where ${\widetilde K}\pm iQ\cos\ff \in C(\pi/2)$.
Setting in \eqref{4.18a}   $K:={\widetilde K}/\sin\ff$
we arrive at  \eqref{4.11}.

It remains to consider the case
$T\in \Ext_{\{T_1,T_2\}}(0)$.
This inclusion means that
$T$ is a self-adjoint contraction in $\gH,$
that is  $T_{21}=T^*_{12}$ and $U=V^*.$
Hence $Q_0=I$ and  $\cH_1=\cH_2=:\cH$.
Threfore equivalence \eqref{4.11}
with $\ff\in(0,\pi/2)$ takes the form
     \begin{equation}\label{4.18b}
T_K\in C_{\gH}(\varphi)\Longleftrightarrow
K\cdot sin\ff\pm i\cos\ff\cdot I\in C_{\cH}(\pi/2).
    \end{equation}
The desired equivalence
      \begin{equation}\label{4.18c}
T_K\in C_{\gH}(0)\Longleftrightarrow K=K^*\in C_{\cH}(\pi/2)
       \end{equation}
is implied now by  \eqref{3.5}.

(iii)  This assertion is immediately implied by
the statement (iii) of Lemma \ref{lem4.2}.
      \end{proof}
       \begin{remark}\label{rem4.3}
Comparison of  condition \eqref{4.10}  with the obvious criterion
$U=V^*$ for the existence of
$T=T^*\in \Ext_{\{T_1,T_1\}}(\pi/2)$
yields a curious fact:
        \begin{equation}\label{4.19A}
U,\ V\in C_\gH(\pi/2),\quad |((I-VU)f,g)|\le
\|D_Uf\|\cdot \|D_{V^*}g\| \equ  U=V^*.
        \end{equation}
I don't know the direct proof if this equivalence.
              \end{remark}

{\bf 3.3. Extreme points of the set $\Ext_{\{T_1,T_2\}}(\ff).$}
Denote by $\Ext^e_{\{T_1,T_2\}}(\ff)$ the set of extreme points of
the closed convex set $\Ext_{\{T_1,T_2\}}(\ff)$.
Theorem \ref{th4.1} makes it
possible to describe a part of the set
$\Ext^e_{\{T_1,T_2\}}(\ff)$.
To this end for any operator
$Q\bigl(\in[\cH_1, \cH_2]\bigr)$ we introduce
the operator loones
   \begin{equation}\label{4.19B}
L(Q;\ff):=\{K\in[\cH_1, \cH_2]:\ K\sin\ff\pm iQ \cos\ff\in C(\pi/2)\},
\qquad \ff\in (0,\pi/2),
   \end{equation}
and denote by $L^e(Q;\ff)$ the set of its extreme points.
     \begin{proposition} \label{prop4.1A}
Let  $Q\in[\cH_1, \cH_2],\  \ff_1 := \arccos(\|Q\|^{-1})>0$
and $\ff\in[\ff_1, \pi/2]$.
Then

(i) the following equivalence holds
   \begin{equation}\label{4.19C}
K\in L(Q;\ff)\equ \sin 2\ff\cdot(K^* Q)_I= D_{K,Q}C D_{K,Q}, \quad
C=C^*\in C_{\cH_1}(0),
   \end{equation}
where $(K^* Q)_I:=(2i)^{-1}(K^* Q-Q^* K)$ and
   \begin{equation}\label{4.19D}
D_{K,Q}:=(I-K^* K\sin^2\ff-Q^* Q\cos^2\ff)^{1/2}\ge 0.
   \end{equation}

(ii) If additionally $\ran D_{K,Q}$ is closed,
that is  $\ran D_{K,Q}= {\overline \ran}D_{K,Q},$
then the following implication holds
    \begin{equation}\label{4.19E}
\sigma(C)\subset \{\pm 1\}\then K\in L^e(Q;\ff).
   \end{equation}
         \end{proposition}
    \begin{proof}
(i) By definition $K\in L(Q;\ff)$ iff
   \begin{equation*}
(I-K^*\sin\ff\mp iQ^*\cos\ff)(K\sin\ff\pm iQ\cos\ff)\ge 0.
   \end{equation*}
With account of definition \eqref{4.19D}
this inequality may be rewritten as
   \begin{equation*}
\pm \sin2\ff(K^* Q)_I\le D^2_{K,Q}.
   \end{equation*}
By Lemma \ref{lem4.1} this inequality is equivalent
to representation \eqref{4.19C} with some selfadjoint
contraction $C$.

(ii) The proof of this statement is similar to that of
Proposition 3.18 from \cite{M4}. Suppose the contrary, that is
$K\notin L^e(Q;\ff).$
Then $2K=K_1+K_2$ where  $K_j\in L(Q;\ff),\  j\in\{1,2\}$
and $K_1\not = K.$
For any $f\in \ker D_{K,Q}$ we have
    \begin{eqnarray*}
4\|f\|^2=\sin^2\ff\cdot\|2K\|^2 +4 \cos^2\ff\cdot\|Q f\|^2
=\sin^2\ff\cdot\|K_1 f+K_2 f\|^2 + 4\cos^2\ff\cdot\|Qf\|^2  \\
\le 2(\sin^2\ff\cdot\|K_1 f\|^2 + \cos^2\ff\cdot\|Q f\|^2)
+ 2(\sin^2\ff\cdot\|K_2 f\|^2 + \cos^2\ff\cdot\|Q f\|^2)
\le 4\|f\|^2.
   \end{eqnarray*}
Hence
   \begin{equation*}
\sin^2\ff\cdot\|K_j f\|^2 + \cos^2\ff\cdot\|Q f\|^2=\|f\|^2,   \qquad
j\in\{1,2\}.
   \end{equation*}
Thus,
$\|K_1 f\|=\|K_2 f\|=\|K f\|$ and $\|K_1 f+K_2 f\|=\|K_1 f\|+\|K_2 f\|.$
In view of strict convexity of the unit ball in $\gH$ we get
    \begin{equation}\label{4.19G}
K_1 f=K_2 f=K f, \qquad   f\in \ker D_{K,Q}.
   \end{equation}

Further, setting
$K_{\pm}:=K \sin\ff\pm iQ\cos\ff$ and using representation
\eqref{4.19C} we obtain
   \begin{equation}\label{4.19K}
D^2_{K_{\pm}}=D^2_{K,Q}\pm \sin 2\ff\cdot(K^* Q)_I=
D_{K,Q}(I\pm C)D_{K,Q}.
   \end{equation}
Suppose that
$f\in(\ker D_{K,Q})^{\perp}$ and $D_{K,Q}f\in \ker(I+C)$.
Then, it follows from \eqref{4.19K}
that $D^2_{K_+}f=0$, that is $\|K_+ f\|=\|f\|$.
Setting
$$
K_{J+}:=K_j\sin\ff + iQ\cos\ff\bigl(\in C(\pi/2)\bigr),\quad  j\in\{1,2\},
$$
and noting that $2K_+=K_{1+}+K_{2+}$ we easily get
    \begin{equation*}
2\|f\|=2\|K_+ f\|=\|(K_{1+}+K_{2+})f\|\le 2\|f\|.
     \end{equation*}
Hence $\|K_{1+}f\|=\|K_{2+}f\|=\|K_+ f\|=\|f\|$ and
   \begin{equation*}
\|K_{1+}f + K_{2+}f\|=\|K_{1+}f\|+\|K_{2+}f\|.
   \end{equation*}
In view of strict convexity of the unit ball in $\gH$ we get
$K_{1+}f=K_{2+}f=K_+ f$, that is
$K_1 f=K_2 f=K f$ for any $f$ obeying $D_{K,Q}f\in \ker(I+Q)$.
Similarly we obtain  that $K_1 f=K_2 f=K f$ for any $f$ such that
$D_{K,Q}f\in \ker(I-Q)$.
Taking into account the hypothesis of proposition we get
   \begin{equation}\label{4.195}
K_1 f=K_2 f=Kf, \quad   f\in\ran D_{K,Q}=ker(I+C)\oplus ker(I-C).
   \end{equation}
Combining \eqref{4.19G} with \eqref{4.195} we get
$K=K_1=K_2.$ This contradicts the assumption that $K_1\not = K.$
    \end{proof}
      \begin{remark}
(a) Closability of the linear manifolds $\ran D_{K,Q}$
in Porposition \ref{prop4.1A}   may be replaced by
${\overline{\ran D_{K,Q}\cap\cH_{\pm}}}=\cH_{\pm}$ where
$\cH_{\pm}:=\ker(I\pm C),$ which are, for example, valid if either
$\dim\cH_+<\infty$ or $\dim\cH_-<\infty$.

(b) Note that $\ran D_{K,Q}$ is closed if both $K$
and $Q$ are compact operators,  $K,Q\in{\gS}_{\infty}$.
      \end{remark}

Next we clarify and complement  Proposition \ref{prop4.1A}
in the case $\cH_1=\cH_2=\cH$  and $Q=I_{\cH}.$ Now
$L(Q;\ff)=L(I_{\cH};\ff)=C_{\cH}(\ff)$. Denote by
$C^e_{\cH}(\ff):=L^e(I_{\cH};\ff)$ the set of extreme points
of the set $C_{\cH}(\ff)$ and by  
     \begin{equation}\label{4.19boundhole}
\partial L_{\ff}:=\partial L^+_{\ff}\cup\partial L^-_{\ff}
\qquad \text{ where}\qquad
\partial L^{\pm}_{\ff}:=\{z\in{\Bbb D}:\  |z\sin\ff\pm i\cos\ff|=1\},
       \end{equation}
the (topological) boundary of the hole \eqref{3.4hole}.
Note that $\partial L_{\ff}$ is at the same time the set
of extreme points of the hole \eqref{3.4hole}.
     \begin{proposition}\label{prop4extreme}
Let $\ff\in(0,\pi/2)$ and $K\in C_{\cH}(\ff)$.
Then

(i) there exists a contraction $C=C^*$ such that
      \begin{equation}\label{4.196}
2 K_I=\tg\ff\cdot D_K C D_K,   \qquad    C\in C_{\cH}(0).
      \end{equation}
Conversly, if $K\in C_{\cH}(\pi/2)$ and \eqref{4.196} holds then
$K\in C_{\cH}(\ff)$;

(ii) the following implication holds
    \begin{equation}
\sigma(C) \subset \{\pm 1\} \quad\text{and}\quad
{\overline{\ran D_{K,Q}\cap\cH_{\pm}}}=\cH_{\pm}:=\ker(I\pm C)
\then K\in C^e_{\cH}(\ff);
     \end{equation}

(iii) if $K$ is a normal operator, $K K^*=K^* K$, and
$\sigma(K)\subset\partial L_{\ff}$
then $K\in C^e_{\cH}(\ff)$;

(iv) if $K\in C_{\cH}(\ff), \  \sigma(K)\subset\partial L_{\ff}$
and the spectrum $\sigma(K)$ is purely point,
then $K$ is normal, hence $K\in C^e_{\cH}(\ff)$.

(v) the set $C^e_{\cH}(\ff)$ contains continuum (nonnormal)
operators  $K$ with $\sigma(K)=0.$
       \end{proposition}
      \begin{proof}
(i) If $Q=I_{\cH}$ then $D_{K,Q}=D_{K,I}=\sin\ff\cdot D_K$ and
the statement is implied by Proposition \ref{prop4.1A} (i).

(ii) This statement is implied by Proposition \ref{prop4.1A} (ii).

(iii) Assume for brevity that
$\pm 1\notin\sigma_p(K)$. Then starting with  \eqref{4.196}
and applying  Spectral theorem we get
   \begin{eqnarray*}
C=\ctg\ff\cdot D^{-1}_K(2 K_I)D^{-1}_K=
\ctg\ff\cdot\int_{\partial L_{\ff}}
\frac{\lambda-{\overline\lambda}}{i\sqrt{1-|\lambda |^2}}d E_K(\lambda)  \\
=\ctg\ff\cdot\int_{\partial L^+_{\ff}}
\frac{\lambda-{\overline\lambda}}{i\sqrt{1-|\lambda |^2}}d E_K(\lambda)+
\ctg\ff\cdot\int_{\partial L^-_{\ff}}
\frac{\lambda-{\overline\lambda}}{i\sqrt{1-|\lambda |^2}}d E_K(\lambda)  \\
= \int_{\partial L^+_{\ff}}d E_K(\lambda)-
\int_{\partial L^-_{\ff}}d E_K(\lambda)=:P_+-P_-.
  \end{eqnarray*}
Here $E_K(\cdot)$ is the spectral measure of $K,$
and  $P_{\pm}$ are the  corresponding   spectral projections.
Since $P_++P_-=I,$ we have
$\sigma(C)\subset \{\pm 1\}.$
Moreover, $P_{\pm}\ran D_K$ is dense in $\cH_{\pm}:=P_{\pm}\cH.$
Hence  by statement (ii) $K\in C^e_{\cH}(\ff)$.

(iv)
Let us set
$K_{\pm}:=K\sin\ff\pm i\cos\ff$.
If $\lambda_j\in\partial L_{\ff}\cap\sigma_p(K)$ and
$\cH_j:=\ker(K-\lambda_j)(\not=\emptyset)$, then either
$\mu^+_j:=\lambda_j\sin\ff+i\cos\ff\in\sigma_p(K_+)\cap {\Bbb T}$ or
$\mu^-_j:=\lambda_j\sin\ff-i\cos\ff\in\sigma_p(K_-)\cap{\Bbb T}$ where
${\Bbb T}:=\{z\in{\Bbb C}:|z|=1\}$.
Thus, the subspace $\cH_j$ reduces the operator $K$
for any $j\in{\Bbb Z}_+$ since either
$\cH_j=\ker(K_+-\mu^+_j)$ or $\cH_j=\ker(K_--\mu^-_j)$ and both
$K_+$ and $K_-$ are contractions. Since the spectrum
$\sigma(K)$ is purely point, then
$K=\oplus^{\infty}_{j=1}\lambda_j I_{\cH_j}$ and $K$ is normal.

(v) First we  consider the case $\cH={\Bbb C}^2$.
We let
$$
K(\theta):=e^{i\theta}
\begin{pmatrix}
0&\sin\ff\\
0&0
\end{pmatrix},
\qquad   \theta\in[0,2\pi].
$$
Then
    \begin{equation*}
D_{K_{(\theta)}}=
\begin{pmatrix}
1&0\\
0&\cos\ff
\end{pmatrix},
\qquad
K(\theta)_I=i\sin\ff
\begin{pmatrix}
0&-e^{i\theta}\\
e^{-1\theta}&0
\end{pmatrix},\qquad
C(\theta)=i
\begin{pmatrix}
0&-e^{-i\theta}\\
e^{-i\theta}&0
\end{pmatrix}.
   \end{equation*}
Hence $\sigma\bigl(C(\theta)\bigr)=\{\pm 1\}$ and by statement (ii)
$K(\theta)\in C^e_{\cH}(\ff)$.
       \end{proof}
     \begin{remark}
(i) Another proof of statement (iii) is contained
in \cite{M4}. The proof of statement (iv) is borrowed
from \cite{M4} and it is presented for the sake if completeness.

(ii) Note that while a complete description of the set
$C^e_{\cH}(\ff)$ is unknown, it essentially differs from
that of the sets  $C^e_{\cH}(0)$ and $C^e_{\cH}(\pi/2)$
even in the case  $\dim\cH<\infty$.
Indeed, if $\dim\cH<\infty$ then by Proposition \ref{prop2.1}
both $C^e_{\cH}(\pi/2)$ and $C^e_{\cH}(0)$ consist of normal
matrices with "boundary spectrum", that is,
$C^e(\pi/2)\ \bigl(\text{resp. }C^e_{\cH}(0)\bigr)$ is
the set of unitary (resp. unitary selfadjoint) matrices.

On the other hand, the sets
$C_{\gH}(\ff), \ff\in(0,\pi)$ may be considered as
"interpolation sets" between $C_{\gH}(0)$ and $C_{\gH}(\pi/2)$.
This observation makes natural the following hypothesis:
for any
$\ff\in(0,\pi)$ the set $C_{\gH}(\ff)$ consists of normal
matrices with "boundary spectrum".
However Proposition   \ref{prop4extreme}
shows that  this hypothesis is
false to be true, since the set $C^e_{\gH}(\ff)$
contains continuum  nonnormal matrices
in addition to the set of normal
matrices with spectrum lying on $\partial L_{\ff}.$
        \end{remark}

Combining Theorem \ref{th4.1} with Proposition  \ref{prop4.1A}
we arrive at the following result.
       \begin{corollary}\label{cor4.1A}
Suppose that conditions of Theorem \ref{th4.1} are satisfied,
$Q_0$ and  $\ff_1$ are  defined by \eqref{4.9} and
\eqref{4.10} respectively, and  $T_K\in \Ext_{\{T_1,T_2\}}(\ff)$.
Then

(i) for any $\ff\in[\ff_1, \pi/2]$ the following equivalence holds
$$
T_K\in \Ext^e_{\{T_1,T_2\}}(\ff)\equ K\in L^e(Q;\ff);
$$

(ii) there exists a selfadjoint contraction $C\in C_{\cH_1}(0)$
such that
$$
\sin 2\ff\cdot(K^* Q)_I=D_{K,Q}C D_{K,Q};
$$

(iii) the following implication holds
$$
\sigma(C)=\{\pm 1\}  \quad \text{and }\quad
\ran D_{K,Q}={\overline\ran}D_{K,Q}\then
T_K\in \Ext^e_{\{T_1,T_2\}}(\ff);
$$

(iv) $T_K\in {\Ext}^e_{\{T_1,T_2\}}(\ff)$
if at least one of the following identities holds
    \begin{equation*}
D^2_{K,Q}\pm\sin 2\ff\cdot(K^*Q)_I=0, \qquad
D^2_{K^*,Q^*}\pm\sin 2\ff(K^*Q)_I=0.
    \end{equation*}
        \end{corollary}

{\bf 3.4. Proper $C_{\gH}(\ff)$-extensions of symmetric
contractions.}
Here we apply Theorem \ref{th4.1} to the case of a dual
pair $\{T_1, T_1\}.$

Let $T_1\in [\gH_1, \gH]$ be a nondensely defined
symmetric contraction in $\gH=\gH_1\oplus\gH_2.$
As usual $\Ext_{T_1}$ stands for the set
of all proper extensions  of $T_1,$ that is
$T\in\Ext_{T_1}$ iff  $T\supset T_1$ and $T^*\supset T_1.$
Denote by
      \begin{equation}\label{4.18A}
\Extp_{T_1}(\ff):=\Ext_{T_1}\cap C_{\gH}(\ff), \qquad  \ff\in [0,\pi/2],
      \end{equation}
the set of all proper $C_{\gH}(\ff)$-extensions
of the symmetric contraction $T_1.$

By Definition \ref{def2.1}  $\Ext_{\{T_1,T_1\}}=\Ext_{T_1}.$
Moreover, it follows from \eqref{4.8A} and \eqref{4.18A}
that  $\Ext_{\{T_1,T_1\}}(\ff)=\Extp_{T_1}(\ff).$
      \begin{corollary}\label{cor4.1}
Let $T_1\in [\gH_1, \gH]$ be a nondensely defined
symmetric contraction in $\gH=\gH_1\oplus\gH_2.$
Then
$\Extp_{T_1}(\ff)=\Ext_{\{T_1,T_1\}}(\ff)\not=\eset$
for any  $\ff\in [0,\pi/2]$. Moreover, the following
equivalence holds
       \begin{equation} \label{4.19}
T:=T_K\in \Extp_{T_1}(\ff)\equ T_{22} = -U^*T_{11}U
+ D_UKD_U,\qquad K\in C_{\cH}(\ff),
        \end{equation}
where $\cH:=\gH_2\ominus \ker D_U$.
     \end{corollary}
         \begin{proof}
According to \eqref{4.8} $T_1=T_2$ if and only if
$T_{12}=T_{21}^*,$  that is iff $V^*=U$.
Therefore the operator $Q_0$ defined by \eqref{4.9}
takes the form $Q_0=I_{\cH}$, where $I_{\cH}$
is the  identical operator in  $\cH.$
Thus $\ff_1=\arccos(\| Q_0\|^{-1})=\arccos 1=0$
and $\Extp_{T_1}(\ff)\not=\eset$ for any $\ff\in [0, \pi/2].$
Moreover, now equivalence
\eqref{4.11} takes the form \eqref{4.19}.
      \end{proof}
       \begin{remark}\label{rem4.2}
In the case $T_1=T_2$ both left and right radii of the balls
$B_1$ and  $B_2$ are equal: $R_l=D_{V^*}=D_U=R_r$.
      \end{remark}

According to \eqref{4.11} the set $\Ext_{T_1}(0)$ of selfadjoint
contractive extensions of $T_1$ forms an operator segment
("the self-adjoint part" of the operator ball
$B(-U^* T_{11}U; D_U,D_U)$) which is parametrized by the
operator segment
$\{K\in[\cH]:-I_{\cH}\le K\le I_{\cH}\}$.

Consider the extremal  selfadjoint contractive extensions
$T_m:=T_{min}$ and $T_M:=T_{max}$ of the operator $T_1.$
It is clear that $T_m:=T_{-I}$ and $T_M:=T_I$  are the extreme points
of the segment  $\Extp_{T_1}(0)$, corresponding to the
operators $K=-I_{\cH}$ and $K=I_{\cH}$ respectively.
Their block-matrix representations are
of the form
      \begin{equation} \label{4.19AA}
T_m=
\begin{pmatrix}
T_{11}&D_{T_{11}}U\\
U^* D_{T_{11}}&-I+U^*(I-T_{11})U
\end{pmatrix}, \qquad
T_M=
\begin{pmatrix}
T_{11}&T_{D_{11}}U\\
U^* D_{T_{11}}&I-U^*(I+T_{11})U
\end{pmatrix}.
   \end{equation}
Using  representations \eqref{4.19AA} we rewrite
description \eqref{4.19} as
      \begin{equation}\label{4.19BB}
T_K\in \Extp_{T_1}(\ff)\equ 2 T_K=(T_M+T_m)+(T_M-T_m)^{1/2}K
(T_M-T_m)^{1/2},\quad K\in C_{\cH}(\ff).
      \end{equation}
Note that this description of the class $\Extp_{T_1}(0)$ has been obtained
by M.G. Krein \cite {Kr1} (see also \cite{AG, KO}).
Other proofs are contained in  \cite{BN2}, \cite{HMS}.
A generalization of the Krein result to the case of
$C_{\gH}(\ff)$-conractions, that is a description of the class
$\Extp_{T_1}(\ff)$ in the form  \eqref{4.19BB} has been
obtained in \cite{AT1, AT2}
(see also \cite{KM1, DM2, M4} for other proofs).

{\bf 3.5. A description of the set of all proper and
improper $C_{\gH}(\ff)$-extensions of  symmetric contractions.}

Let  $A$  be  a closed densely defined  symmetric
operator in $\gH.$
It is known, that any $m$-dissipative (in particular selfadjoint)
extension  ${\wt A}$ of $A$ is a proper extension
$({\wt A}\in \Extp_A),$ that is
$A\subset {\wt A}\subset A^*$.
It is not the case
for $m$-sectorial extensions of
a nonnegative operator $A\ge 0.$

Therefore we clarify Definition \ref{def3.2A}
for the case of a nonnegative operator.
     \begin{definition}\label{def3.2B}
Let  $A(\ge 0)$  be  a closed densely defined  nonnegative
operator in $\gH.$
Denote by $\Ext_A((0,\infty);\ff)$  the class
of all proper $m$-sectorial
extensions of  $A$ with  vertex zero and half-angle
$\ff\in (0;\pi/2).$
The class of all (proper and improper) $m$-sectorial
extensions ${\wt A}(\in \cC(\gH))$ of $A$  will be
denoted by  ${\Extp}_A((0,\infty);\ff).$
      \end{definition}

Here we present a description of the set
${\Ext}_A((0,\infty);\ff).$
In accordance with the approach accepted in this paper
(cf. Lemma \ref{lem3.2}) it  suffices to describe
the set ${\Ext}_{T_1}(\ff)$  of all the extensions
of the class $C_{T_1}(\ff)$
of a nondensely defined Hermitian contraction
$T_1:=(I-A)(I+A)^{-1}(\in [\gH_1,\gH])$
where $\gH_1:=\ran(I-A).$
In turn,
considering the block-matrix
representation of $T_1$ with respect to the orthogonal
decomposition $\gH=\gH_1\oplus\gH_2,$
one reduces the problem to the problem
of  a description of all the "completions"
of a contractive operator-matrix
$T_1=\binom{T_{11}}{T_{21}}$ to form a matrix
$T=(T_{ij})_{i,j=1}^2\in C_{\gH}(\ff)$.

Note that ${\wt A}\in \Extp_A$ iff the entries
$T_{12}$ and $T_{21}$ of
$(T_{i j})^2_{i,j=1}:=T:=(I-{\wt A})(I+{\wt A})^{-1}$ are
connected by $T_{21}=T^*_{12}$.

       \begin{theorem}\label{th4.2}
Let $T_1=\binom{T_{11}}{T_{21}}=
\binom{T_{11}^*}{VD_{T_{11}}}$  be a symmetric
contraction in $\gH=\gH_1\oplus \gH_2,\ $
$T:=(T_{ij})_{i,j=1}^2\in [\gH],$
$\cH_2:={\overline\ran}D_{V^*},\ $
$\cH'_2:={\overline\ran}D_V,$
and $\ff\in[0,\pi/2].$
Then

(i) $T\in{\Ext}_{T_1}(\pi/2),$  i.e.  $T(\in C_{\gH}(\pi/2))$
is a contractive extension of $T_1$  if and only if
it is of the form \eqref{2.15}, that is
      \begin{equation} \label{4.20}
T_{12}=D_{T_{11}}U,\quad T_{22}=-VT_{11}^*U+D_{V^*}KD_U\quad
\text{with} \quad U,K\in C(\pi/2).
      \end{equation}

(ii) $T\in {\Ext}_{T_1}(\ff):=
{\Ext}_{T_1}(\pi/2)\cap C_\gH(\ff)$
if and only if $\  U$
"runs through" the operator ball of the form
     \begin{equation} \label{4.21}
\begin{split}
U=\sin\ff&(\sin^2\ff D_V^2+\cos^2\ff)^{-1/2}D_VMD_{V^*}(\sin^2\ff D_{V^*}^2
+\cos^2\ff)^{-1/2}+\\
&+\cos^2\ff(\sin^2\ff D_V^2+\cos^2\ff)^{-1}V^*,\qquad M\in C(\pi/2)\cap
[\cH_2, \cH'_2],
\end{split}
      \end{equation}
and $T_{22}$ (for fixed $U$) "runs through" the operator "hole"
       \begin{equation} \label{4.22}
T_{22}=-VT^*U+D_{V^*}KD_U,\qquad K\sin\ff\pm iQ\cos\ff \in C(\pi/2)\cap
[\cH_1,\cH_2].
        \end{equation}
Here $\cH_1:=\cH_1(U):={\overline\ran}D_{U},\ $
$Q={\overline Q_0}(\in [\cH_1,\cH_2])$
and $Q_0$ is defined by  \eqref{4.9}.
      \end{theorem}
    \begin{proof}
Equality \eqref{4.20} is implied by  Theorem \ref{th2.1}.

Let further $T_{12}=D_{T_{11}}U,\ T_2=\binom{T_{11}}{U^*D_{T_{11}}}$.
Then $\dup$ is a dual pair of contractions and according
to Theorem \ref{th4.1} the condition
$\Ext_{\{T_1,T_2\}}(\ff) \not=\eset$
is equivalent to the contractibility of the
operator $Q_0\cos\ff,$  where $Q_0$ is defined by
\eqref{4.9}, i.e.  to the inequality
       \begin{equation} \label{4.23}
\cos\ff\cdot \|D_{V^*}^{-1}(I-VU)f\|\le \|D_Uf\|,\qquad  f\in \gH_2.
        \end{equation}
Supposing first that $0\in\rho(D_{V^*})$, we rewrite  inequality
\eqref{4.23} in the equivalent form
         \begin{equation*}
\cos^2\ff[D_{V^*}^{-2}-U^*V^*D_{V^*}^{-2}-
D_{V^*}^{-2}VU+U^*V^*D_{V^*}^{-2}VU]\le I-U^*U,
          \end{equation*}
or
          \begin{equation} \label{4.24}
     \begin{split}
U^*(1+\cos^2\ff V^*D_{V^*}^{-2}V)U-\cos^2\ff&(U^* D_V^{-2}V^*+VD_V^{-2}U)\\
&+\cos^2\ff D_{V^*}^{-2}-I\le 0.
     \end{split}
           \end{equation}
Since inequality \eqref{4.24} is equivalent to \eqref{4.23}
the set of its solutions is nonempty for any fixed $V.$
By Lemma \ref{lem2.2} for any fixed $V$
the set of solutions of inequality
\eqref{4.24}, that is the set of operators $U$ obeying \eqref{4.24},
forms an operator ball $B(C_0;R_l,R_r)$.
Applying  Lemma \ref{lem2.2} we find
its center and radii. We have
       \begin{equation} \label{4.25}
    \begin{split}
      R_l
        &= Q_1^{-1/2}
        = (I+\cos^2\ff V^*D_{V^*}^{-2}V)^{-1/2}
        = (I+\cos^2\ff V^*VD_{V}^{-2})^{-1/2} \\
        &= [\cos^2\ff(D_V^{-2}-I)+I]^{-1/2}
        = [D_V^{-2}(\sin^2\ff D_V^2+\cos^2\ff)]^{-1/2} \\
        &= (\sin^2\ff D_V^2+\cos^2\ff)^{-1/2}D_V,
    \end{split}
  \end{equation}
     \begin{equation}\label{4.26}
\begin{split}
      C_0
        &= -Q_1^{-1}Q_2
        = -R_l^2Q_2
        = \cos^2\ff(\sin^2\ff D_V^2+\cos^2\ff)^{-1}
          D_V^2D_V^{-2}V^* \\
        &= \cos^2\ff V^*(\sin^2\ff D_{V^*}^2+\cos^2\ff)^{-1}.
    \end{split}
\end{equation}
And finally
       \begin{equation*}
  \begin{split}
    R_r^2
      &= Q_2^*Q_1^{-1}Q_2-Q_3
      = -Q_2^*C_0-Q_3 \\
      &= \cos^4\ff VD_V^{-2}V^*
        (\sin^2\ff D_{V^*}^2+\cos^2\ff)^{-1}
        I-\cos^2\ff D_{V^*}^{-2} \\
      &= \cos^4\ff(I-D_{V^*}^2)D_{V^*}^{-2}
        (\sin^2\ff D_{V^*}^2+\cos^2\ff)^{-1}
        +I-\cos^2\ff D_{V^*}^{-2} \\
      &= (\sin^2\ff D_{V^*}^2+\cos^2\ff)^{-1}\cdot \\
        &\qquad
       \cdot [\cos^4\ff D_{V^*}^{-2}-\cos^4\ff\cdot I
          +(I-\cos^2\ff D_{V^*}^{-2})(\sin^2\ff D_{V^*}^2
          +\cos^2\ff)] \\
      &= (\sin^2\ff D_{V^*}^2+\cos^2\ff)^{-1}
        [\sin^2\ff D_{V^*}^2
          +\cos^2\ff(1-\cos^2\ff-\sin^2\ff)] \\
      &= \sin^2\ff
        (\sin^2\ff D_{V^*}+\cos^2\ff\cdot I)^{-1}D_{V^*}^2.
  \end{split}
\end{equation*}
Thus,
   \begin{equation} \label{4.27}
R_r=\sin\ff(\sin^2\ff D_{V^*} + \cos^2\ff\cdot I)^{-1/2}D_{V^*}.
       \end{equation}
Applying Lemma \ref{lem2.2} and taking relations \eqref{4.25}-\eqref{4.27}
into account we get that inequality \eqref{4.24} or,
what is the same, inequality \eqref{4.23}
is satisfied iff $U$ admits
a representation \eqref{4.21} with some $M\in C(\pi/2)$.
Thus, we proved \eqref{4.21} under the additional assumption
$0\in\rho(D_{V^*})\bigl(\equ 0\in\rho(D_V)\bigr)$.


Next, we may easily free ourselves of the additional
assumption $0\in \rho(D_{V^*})$ by passing to the limit.
Actually since $D_{V^*}^{-1}>
D_{rV^*}^{-1},\quad r\in (0,1)$,
inequality \eqref{4.23} takes place if and only if for any
 $r\in (0,1)$
the inequality
              \begin{equation} \label{4.28}
\cos\ff\|D_{rV^*}^{-1}(I-V)f\|\le \|D_Uf\|,\qquad f\in \gH_2,
               \end{equation}
holds true.
Since $0\in \rho(D_{rV})$, then in accordance with what has been proved
in the previous step
inequality \eqref{4.28} (for fixed $r<1$)
is equivalent to  equality \eqref{4.21} with
$D_V$ and $D_{V^*}$ repleaced by $D_{rV}$ and $D_{rV^*}$
respectively.
In these equalities it is possible to pass to the limit
as $r\uparrow 1$ (in the sence of strong convergence).

Now the relations \eqref{4.22} follow from Theorem \ref{th4.1}.
       \end{proof}
     According to Theorem \ref{th4.2} (see formulas
\eqref{4.20}-\eqref{4.22}) any extension
$T\in{\Ext}_{T_1}(\ff)$ is uniquely determined by a pair
$\{M,K\}$ of "free" parameters. Denote the corresponding extension
$T$ by $T_{M,K}$.

Next we denote by ${\Ext^e}_{T_1}(\ff)$ the set of extreme points
of ${\Ext}_{T_1}(\ff)$.
 \begin{corollary}
Let $T=T_{M,K}\in{\Ext}_{T_1}(\ff)$.
Then

(i) $T_{M,K}\in{\Ext}^e_{T_1}(\ff)$
if and only if $M$ is a maximal partial isometry
from $\cH_2$ to $\cH'_2$ and $K\in L^e(Q;\ff)$,
where $L(Q;\ff)$ is defined by  \eqref{4.19B};

(ii) if $M$ is a maximal partial isometry, then
the following implication holds
$$
\sigma(C)=\{\pm 1\}  \quad \text{and }\quad
\ran D_{K,Q}={\overline\ran}D_{K,Q}\then
T_{M,K}\in {\Ext}^e_{\{T_1,T_2\}}(\ff).
$$
  \end{corollary}
   \begin{proof}
It follows from \eqref{4.20} and \eqref{4.21} that the mapping
$\{M,K\}\to T_{M,K}$ preserves convexity:
if the  "free" parameters
$\{M_j,K_j\},\  j\in\{1,2\}$ and $\{M,K\}$ are connected by
$M=tM_1 + (1-t)M_2$ and $K= tK_1 + (1-t)K_2$ with $t\in (0,1),$
then  $T_{M,K}=tT_{M_1,K_1} + (1-t)T_{M_2,K_2}.$
 To complete the proof it remains to apply both
Proposition \ref{prop2.1}  and Corollary \ref{cor4.1A}.
   \end{proof}
      \begin{remark}   
(i)
The set ${\Ext}_{T_1}(\ff)$ of all (proper and improper)
$C_{\gH}(\ff)$-extensions of $T_1$ admits a representation
${\Ext}_{T_1}(\ff)=\cup_{T_2}\Ext_{\{T_1,T_2\}}(\ff)$
where
$T_2=
\binom {T_{11}}
{U^*D_{T_{11}}}$
and $U$ "runs through" the operator ball \eqref{4.21}.
Note that $T$ is a proper $C_{\gH}(\ff)$-extension of $T_1$ iff
$U=V^*$. In this case $Q=I$ and \eqref{4.22} turns into
\eqref{4.19}.

(ii) Theorem \ref{th4.2} has been proved by the author
together with V. Kolmanovich in \cite{KM1}
in a different but equivalent form.
     \end{remark}

{\bf 3.6.  $C_\gH(\ff)$-extensions
of a dual pair of  $C(\ff)$-contractions.}
Here we consider a dual pair $\{T_1,T_2\}$ of
$C(\ff)$-contractions of the form
    \begin{equation}\label{2.14new}
T_1=\binom{T_{11}}{T_{21}}=\binom{T_{11}}{VD_{T_{11}}},\quad
T_2=\binom{T_{11}^*}{T_{12}}=\binom{T_{11}^*}{U^*D_{T_{11}^*}},
     \end{equation}
and show that the Problem $2C$ mentioned in the
Introduction is reduced to the Problem 3 with
different left $R^+_l\not =R^-_l$  and right
$R^+_r\not =R^-_r$ radii.
        \begin{proposition}\label{prproblem4}
Let $\{T_1,T_2\}$ be a dual pair of conractions in
$\gH=\gH_1\oplus\gH_2$.
Suppose that both $T_1$ and $T_2$ obey condition  \eqref{3.3}
with $\ff=\ff_0.$  Then

(i) for any $\ff\in[\ff_0,\pi/2)$ the following relations hold
      \begin{equation}\label{4.28A}
2(T^*_{11})_I=\tg\ff\cdot D_{T_{11}}C(\ff)D_{T_{11}}=
\tg\ff\cdot D_{T^*_{11}}C_2(\ff)D_{T_{11}^*},
      \end{equation}
where $C_j(\ff):=C_j\cdot\tg\ff_0/\tg\ff,\  j\in\{1,2\},$
is a selfadjoint  conraction  and
       \begin{equation}\label{4.28B}
C_1=D_V C' D_V, \quad   C_2=D_{U^*}C'' D_{U^*}, \quad
 -I\le C',\ C'' \le +I;
       \end{equation}

(ii) for any $\ff\in(\ff_0,\pi/2)$ the set
$\Ext_{\{T_1,T_2\}}(\ff)$ forms
an operator hole:
      \begin{equation}
T\in\Ext_{\{T_1,T_2\}}(\ff)\equ T_{22}\sin\ff\in L(\ff):=
B(C_+;  R^+_l,R^+_r)\cap B(C_-; R^-_l,R^-_r),
    \end{equation}
where
     \begin{equation}
C_{\pm}:=\mp i\cos\ff\cdot I - V_{\pm}(T^*_{11}\sin\ff\pm i\cos\ff)U_{\pm},
\quad  R^{\pm}_l=D_{V^*_{\pm}},\  R^{\pm}_r=D_{U_{\pm}}
    \end{equation}
and $V_{\pm},\  U_{\pm}$ are contractions of the form
     \begin{equation}\label{4.28+}
V^*_{\pm}:=V^*_{1\pm}\bigl(I\pm C_1(\ff)\bigr)^{-1/2}V^*, \qquad
U_{\pm}:=U_{1\pm}\bigl(I\pm C_2(\ff)\bigr)^{-1/2}U,
      \end{equation}
and $V_{1\pm},\  U_{1\pm}$ are (uniquely determined) partial isometries.

In particular, $\Ext_{\{T_1,T_2\}}(\ff)\not =\emptyset$ if and only if
$L(\ff)\not = \emptyset.$
        \end{proposition}
   \begin{proof}
The inclusion $T=(T_{ij})^2_{i,j=1}\in C_{\gH}(\ff)$ means that
  \begin{equation}\label{4.28AA}
    T_{\pm} :=
    \begin{pmatrix}
      B_{11}^{\pm} & B_{12} \\
      B_{21}       & B_{2}^{\pm}
    \end{pmatrix}
    :=
    \begin{pmatrix}
      \sin\ff\cdot T_{11}\pm i\cos\ff\cdot I & T_{21}\sin\ff \\
      T_{21}\sin\ff & \sin\ff\cdot T_{22}\pm i\cos\ff\cdot I
    \end{pmatrix}
    \in C_{\gH}(\pi/2).
  \end{equation}
First we note that
   \begin{equation*}
D^2_{T_1}=I-T^*_{11}T_{11}-T^*_{21}T_{21}=
D^2_{T_{11}}-D_{T_{11}}V^* V D_{T_{11}}=
D_{T_{11}}D^2_V D_{T_{11}}
   \end{equation*}
and $D^2_{T_2}=D_{T^*_{11}}D^2_{U^*}D_{T^*_{11}}$.
Combining these relations
with \eqref{3.3} and applying
Lemma \ref{lem4.1} we obtain \eqref{4.28A}.

Next, starting with \eqref{4.28AA}
and taking \eqref{4.28A} into account we get
      \begin{eqnarray}\label{4.28B}
D^2_{B^{\pm}_{11}}=I-(\sin\ff\cdot T^*_{11}\mp i\cos\ff)
(\sin\ff\cdot T_{11}\pm i\cos\ff)  \\
\nonumber
= \sin^2\ff\cdot D^2_{T_{11}}\mp\sin 2\ff\cdot(T_{11})_I=
\sin^2\ff\cdot D_{T_{11}}\bigl(I\pm C_1(\ff)\bigr)D_{T_{11}}.
      \end{eqnarray}
Hence there exist partial isometries $V_{1\pm}(\ff)$ such that
       \begin{equation}\label{4.28C}
D_{B^{\pm}_{11}}=\sin\ff\cdot D_{T_{11}}\bigl(I\pm C_1(\ff)\bigr)^{1/2}
V_{1\pm}(\ff)=\sin\ff\cdot V^*_{1\pm}(\ff)(I\pm C_1)^{1/2}D_{T_{11}}
      \end{equation}
Combining \eqref{4.28AA}, \eqref{4.28B} and
\eqref{4.28C} we derive
     \begin{equation*}
B_{21}=T_{21}\sin\ff=\sin\ff\cdot V D_{T_{11}}=V_{\pm}D_{B^{\pm}_{11}}=
\sin\ff\cdot V_{\pm}V^*_{1\pm}(\ff)(I\pm C_1)^{1/2}\cdot  D_{T_{11}}
      \end{equation*}
It follows that
$V=V_{\pm}\cdot V^*_{1\pm}(\ff)\bigl(I\pm C_1(\ff)\bigr)^{1/2},$
which yields the first of relations \eqref{4.28+}.

Similarly we get
      \begin{equation}
D^2_{B^{\pm *}_{11}}=\sin^2\ff\cdot
D_{T^*_{11}}\left(I\pm C_2(\ff)\right)D_{T^*_{11}}.
    \end{equation}
According to polar decomposition we have
$\sin\ff\cdot\left(I\pm C_2(\ff)\right)^{1/2}D_{T^*_{11}}=
U_{1\pm}D_{B^{\pm *}_{11}}$ with some partial isometries $U_{\pm}$.
These representations imply
      \begin{equation}
B_{12}=T_{12}\sin\ff=\sin\ff\cdot D_{T^*_{11}}U=
D_{B^{\pm *}_{11}}U_{\pm}=
\sin\ff\cdot D_{T^*_{11}}\left(I\pm C_2(\ff)\right)^{1/2}U_{1\pm}U_{\pm}.
    \end{equation}
Hence $U=\left(I\pm C_2(\ff)\right)^{1/2}U_{1\pm}U_{\pm}$.
This equality yields the second relation in \eqref{4.28+}.

By Theorem \ref{2.4}  $T_{\pm}\in C_{\gH}(\pi/2)$ if and only if
      \begin{equation}
T_{22}\sin\ff\pm i\cos\ff\cdot I = C'_{\pm}+D_{V^*_{\pm}}K_{\pm}D_{U_{\pm}},
      \end{equation}
where $K_{\pm}$ are contractions and $C'_{\pm}$ are defined by
        \begin{equation}
C'_{\pm}= -V_{\pm}(B^{\pm}_{11})^* U_{\pm}=
-V_{\pm}(T^*_{11}\sin\ff\pm i\cos\ff\cdot I)U_{\pm}.
        \end{equation}
to complete the proof it suffices to
set $C_{\pm}=C'_{\pm}\mp i\cos\ff\cdot I.$
     \end{proof}
        \begin{remark}\label{remProblem3}
Let $T_1=\binom{T_{11}}{T_{21}}\in C(\ff)$ and
$\gH=\gH_1\oplus\gH_2$.
It is shown in \cite{M4}, Theorem 4.11,
that Problem 3C is reduced to Problem 5 mentioned
in the Introduction. Namely, it is proved in \cite{M4} that
$T\in\Ext_{T_1}(\ff)$ iff
   \begin{equation}
T P_2\in L:=B_+\cap B_-,  \qquad \text{where} \qquad
B_{\pm}=B(C_{\pm}; D_{S^*_{\pm}}/\sin\ff, P_2),
   \end{equation}
and $S_{\pm}=T_1\sin\ff\pm i\cos\ff\cdot I,\ \
 C_{\pm}=\mp\ctg\ff\cdot P_2.$
Thus, the set  $\Ext_{T_1}(\ff)$ forms an operator
hole of the form \eqref{1.3} with $R^{\pm}_l=D_{S^*_{\pm}}/\sin\ff,$
$R^{\pm}_r=P_2,$ and $C_{\pm}=\mp\ctg\ff\cdot P_2.$
      \end{remark}

\section{Noncontractive extensions of dual pair of symmetric
contractions.}

{\bf 4.1. Schur complements.}
In this section we investigate some spectral properties of
contractive and noncontractive extensions
of a dual pair $\{T_1, T_2\}$ of symmetric  contractions
using their block-matrix representations \eqref{4.8}.
Trough this section we keep a notation $T_K$
for any (not necessary contractive)
extension of the dual pair $\{T_1, T_2\}$ having the form
\eqref{2.15} with a bounded operator $K\in[{\cH}_1, {\cH}_2].$
Observe  that any bounded extension $T\in \Ext_{\{T_1,T_2\}}$
has such a form iff  $T_1$ and   $T_2$  are transversal,
that is  $0\in \rho(D_U)\cap \rho(D_{V^*}).$
Note also that in the nonsingular
case  $\cH_1=\gH_1$ and $\cH_2=\gH_2'.$

We investigate some spectral properties of extensions
$T_K(\in \Ext_{\{T_1,T_2\}})$
in terms of "boundary" operators $K.$
In particular we obtain  descriptions of the
classes $C_{\gH}(\ff;\kk^{\pm})$  and $C_{\gH}(\ff;\gS^{\pm}).$

As well as in Theorem \ref{th4.2} these descriptions
essentially depend on the operator
    \begin{equation} \label{4.9A}
Q_0:= D_{V^*}^{-1}(I-VU)D_U^{-1}.
     \end{equation}

In the following theorem which  is the main result of the
section we calculate Schur complement of the operator
block-matrices   $I-T_KT_K^*\pm \ctg\ff\cdot(T_K-T_K^*).$
     \begin{theorem}\label{th4.3}
Let $\dup$ be a dual pair of symmetric  contractions
of the form \eqref{4.8} in $\gH=\gH_1\oplus\gH_2,\ $
$T_K\in \Ext_{\{T_1,T_2\}}$ and let
$\ff\in [\ff_1,\pi/2]$  where $\ff_1:=\arccos({\|Q\|}^{-1}),\ $
 $Q:={\overline Q_0}$ and  $Q_0$ be defined by \eqref{4.9A}.
Let further
           \begin{equation} \label{4.29}
S^{\pm} := (S^{\pm}_{ij})_{i,j=1}^2 :=I-T_KT_K^*\pm \ctg\ff\cdot(T_K-T_K^*),
      \end{equation}
          \begin{equation}\label{4.30}
G^{\pm} := (G^{\pm}_{ij})_{i,j=1}^2 :=I-T_K^*T_K\pm \ctg\ff\cdot(T_K-T_K^*),
           \end{equation}
with $S_{ij}^{\pm},\ G_{ij}^{\pm}\in [\gH_i,\gH_j]$. Then

(i) $\ran(S_{11}^{1/2})\supset \ran(S_{12}^{\pm}),\
\ran(G_{11}^{1/2})\supset \ran(G_{12}^{\pm})$ and,
consequently the operators
$S_{11}^{-1/2}S_{12}^{\pm}\ $ and \  $G_{11}^{-1/2}G_{12}^{\pm}$
are well  defined and bounded, where $S_{11}:=S_{11}^+=S_{11}^-;$

  (ii)
  the identities
      \begin{eqnarray} \label{4.31}
    &\sin^2\ff&\cdot[S_{22}^{\pm}-(S_{11}^{-1/2}S_{12}^{\pm})
            (S_{11}^{-1/2}S_{12}^{\pm})^*] \\
       \nonumber
 & = & D_{V^*}\cdot[I-(K\sin\ff\mp iQ\cos\ff)
(K^*\sin\ff\pm i Q^*\cos\ff)]\cdot D_{V^*},
      \end{eqnarray}
      \begin{eqnarray}    \label{4.32}
   &\sin^2\ff&\cdot[G_{22}^{\pm}-(G_{11}^{-1/2}G_{12}^{\pm})^*
            (G_{11}^{-1/2}G_{12}^{\pm})] \\
      \nonumber   
  & = & D_U\cdot [I-(K^*\sin\ff\pm iQ^*\cos\ff)
(K\sin\ff\mp iQ\cos\ff)]\cdot D_U,
  \end{eqnarray}
  hold true.
          \end{theorem}
        \begin{proof}
(i) We let
$(G_{ij})_{i,j=1}^2:=I-T_K^*T_K$. Then
        \begin{equation} \label{4.33}
(G_{ij}^{\pm})_{i,j=1}^2=(G_{ij})_{i,j=1}^2\pm i\ctg \ff
\begin{pmatrix}
0&D_{T_{11}}(U-V^*)\\
(V-U^*)D_{T_{11}}&T_{22}-T_{22}^*
\end{pmatrix}.
        \end{equation}
By definition the operator
 $T_K (\in \Ext_{\{T_1,T_2\}})$ is of the form \eqref{2.15}
with  $T_{22}=-VT_{11}U+D_{V^*}KD_U$ and  $K\in[{\cH}_1, {\cH}_2].$
Therefore taking into account \eqref{2.15} and
\eqref{4.33} we get
      \begin{equation} \label{4.34}
G_{11}^{\pm}=G_{11}=D_{T_{11}}D_V^2D_{T_{11}},
      \end{equation}
and
        \begin{equation} \label{4.35}
-G_{21}^{\pm}=-(G_{12})^{\pm})^*=
(U^*T_{11}D_V+D_UK^*V)D_VD_{T_{11}}\mp i\ctg \ff\cdot(U^*-V)D_{T_{11}}.
      \end{equation}
Since $\ff\in [\ff_1, \pi/2]$, then according to Theorem \ref{th4.2}
the operators $U$ and $V$ are connected by
equality \eqref{4.21}.  Setting
      \begin{equation} \label{4.36}
Y:=(\sin^2\ff D_V^2+\cos^2\ff)^{-1/2} \quad \text{and}\quad
Y_*=(\sin^2\ff D_{V^*}^2 + \cos^2\ff)^{-1/2},
       \end{equation}
and taking into account the identity  $VD_V=D_{V^*}V$
we rewrite \eqref{4.21} in the form
        \begin{equation} \label{4.37}
U^*-V=\sin\ff D_{V^*}Y_*\cdot\bigl(M^*-\sin\ff\cdot V\bigr)\cdot YD_V.
        \end{equation}
Let, further
       \begin{equation} \label{4.38}
X_{\pm} := U^*T_{11}D_V + D_UK^*V\mp i\cos\ff\cdot D_{V^*}Y_*\cdot
\bigl(M^*-\sin\ff\cdot V\bigr)\cdot Y.
       \end{equation}
Now  relations \eqref{4.35}-\eqref{4.38} yield
$-G_{21}^{\pm}=X_{\pm}D_VD_{T_{11}}.$  Combining this equality
with  \eqref{4.34} we easily get
     \begin{equation} \label{4.39}
\|G_{12}^{\pm}f\|^2 \le \|X_{\pm}\|^2\cdot (G_{11}f,f)=
 \|X_{\pm}\|^2 \cdot \|(G_{11}^{1/2}f\|^2,
\qquad  f\in \gH_1.
      \end{equation}
This inequality yields the inclusion
$\ran(G_{12}^{\pm})\subset \ran(G_{11}^{1/2}),$
that is the second  of the required inclusions.

The proof of the first  inclusion
$\ran(S_{11}^{12})\supset \ran(S_{12}^{\pm})$
can be obtained in just the same way.
It is suffices to use the equalities
$S_{11}^{\pm}=S_{11}=D_{T_{11}}D_{U^*}^2D_{T_{11}}$
and $S_{21}^{\pm}=X_{\pm}'D_{U^*}D_{T_{11}}$
in place of  \eqref{4.34}  and \eqref{4.38} respectively,
and the relation
         \begin{multline} \label{4.37'}
  V-U^*
  = \sin\ff \cdot D_U(\sin^2\ff D_U^2+\cos^2\ff)^{-1/2}\cdot  \\
    \bigl(M_1-\sin\ff\cdot U^*\bigr)\cdot
 (\sin^2\ff D_{U^*}^2+\cos^2\ff)^{-1/2}D_{U^*}
         \end{multline}
in place of \eqref{4.37}.

$(ii_1)$  Let us prove equality \eqref{4.32}
assuming at the begining that
$0\in\rho(D_{T_{11}})\cap\rho(D_V)$.
In this case setting  $Z:=T_{22}$, we obtain
from \eqref{4.33} - \eqref{4.35} that
    \begin{equation} \label{4.40}
\begin{split}
&G_{22}^{\pm}-(G_{12}^{\pm})^*G_{11}^{-1}G_{12}^{\pm} =
I-U^*D_{T_{11}}^2U-Z^*Z\pm i\ctg\ff\cdot (Z-Z^*)\\
&-[U^*(T_{11}\pm i\ctg\ff)\pm (Z^*\mp i\ctg\ff)V]\cdot
D_V^{-2} \cdot [(T_{11}\mp i\ctg \ff)U+V^*(Z\pm i\ctg\ff)]\\
&=(1+\ctg^2\ff)\cdot I-U^*D_{T_{11}}^2U-(Z^*\mp i\ctg\ff\cdot I)\cdot
(Z\pm i\ctg\ff\cdot I)\\
&-(Z^*\mp i\ctg \ff)\cdot VD_V^{-2}V^*\cdot (Z\pm i\ctg \ff)
- U^*(T_{11}\pm i\ctg\ff) \cdot D_V^{-2}V^* \cdot (Z\pm i\ctg\ff)\\
&-(Z^*\mp i\ctg\ff)\cdot VD_V^{-2} \cdot(T_{11}\mp i\ctg\ff)U-
U^*(T_{11}\pm i\ctg\ff)\cdot D_V^{-2}\cdot (T_{11}\mp i\ctg\ff)U\\
&=D_U^2+\ctg^2\ff\cdot I+U^*T_{11}^2U-(Z^*\mp i\ctg\ff)\cdot D_{V^*}^{-2}\cdot
(Z\pm i\ctg\ff\cdot I)\\
&-U^*(T_{11}\pm i\ctg\ff)\cdot D_V^{-2}V^* \cdot (Z\pm i\ctg\ff)-
(Z^*\mp i\ctg\ff)\cdot VD_V^{-2}\cdot (T_{11}\mp i\ctg\ff)U \\
&-U^*(T_{11}\pm i\ctg\ff)\cdot D_V^{-2} \cdot (T_{11}\mp i\ctg\ff)U.
\end{split}
     \end{equation}

On the other hand, combining \eqref{4.9A} with
the equality $Z:=T_{22}=-VT_{11}U+D_{V^*}KD_U$, we get
    \begin{equation} \label{4.41}
QD_U=D_{V^*}^{-1}(I-VU) \qquad \text{and} \qquad
KD_U=D_{V^*}^{-1}(Z+VT_{11}U).
     \end{equation}
Inserting these relations in the right-hand side of
\eqref{4.32} we deduce
     \begin{equation} \label{4.42}
\begin{split}
&A_{\mp} :=D_U\cdot \left[\frac{1}{\sin^2\ff}-
\left(K^*\pm i{Q^*}{\ctg\ff}\right)\cdot
\left(K\mp iQ{\ctg\ff}\right)\right]\cdot D_U  \\
&= \frac {D_U^2}{\sin^2\ff}
-[(Z^*+U^*T_{11}V^*)\mp i\ctg\ff(I-U^*V^*)]\cdot
D_{V^*}^{-2}\cdot [(Z+VT_{11}U)\pm i\ctg\ff(I-VU)]\\
&=\frac{D_U^2}{\sin^{2}\ff}-[(Z^*\mp i\ctg\ff)+U^*(T_{11}\pm i\ctg\ff)V^*]\cdot
D_{V^*}^{-2} \cdot [(Z\pm i\ctg\ff)+V(T_{11}\mp i\ctg\ff)U]\\
&=\frac{D_U^2}{\sin^2\ff}-(Z^*\mp i\ctg\ff)\cdot D_{V^*}^{-2}\cdot
(Z\pm i\ctg\ff)-(Z^*\mp i\ctg\ff)\cdot D_{V^*}^{-2}V \cdot
(T_{11}\mp i\ctg\ff)U\\
&-U^*(T_{11}\pm i\ctg\ff)\cdot V^*D_{V^*}^{-2}\cdot
(Z\pm i\ctg\ff)-U^*(T_{11}\pm i\ctg\ff)
\cdot V^*D_{V^*}^{-2}V \cdot (T_{11}\mp i\ctg\ff)U.
\end{split}
     \end{equation}
Since $V^*D_{V^*}^{-2}V=D_V^{-2}V^*V=D_V^{-2}-I$,
then the last term in \eqref{4.42} is transformed as follows:
     \begin{equation} \label{4.43}
\begin{split}
-U^*(T_{11}\pm i\ctg\ff)\cdot V^*D_{V^*}^{-2}V\cdot (T_{11}\mp i\ctg\ff)U
=U^*T_{11}^2U\\
+\ctg^2\ff\cdot U^*U-U^*(T_{11}-i\ctg\ff) \cdot D_V^{-2}\cdot
(T_{11}+i\ctg\ff)U.
\end{split}
     \end{equation}
Comparing \eqref{4.40} with \eqref{4.42} and \eqref{4.43} and noting that
      \begin{equation*}
\frac1{\sin^2\ff}D_U^2+\ctg^2\ff\cdot U^*U=
D_U^2+\ctg^2\ff \cdot (D_U^2+U^*U) = D_U^2+\ctg^2\ff\cdot I,
      \end{equation*}
we arrive at the equality
      \begin{equation*}
A_{\mp}=G_{22}^{\pm}-(G_{11}^{-1/2}G_{12}^{\pm})^*(G_{11}^{-1/2}G_{12}^{\pm}),
       \end{equation*}
coinciding with \eqref{4.32}.

$(ii_2)$ Now we free ourselves of the additional restriction
$0\in \rho(D_{V^*})\cap \rho(D_{T_{11}})$.
Consider the strict contractions
$rT_1= \binom{rT_{11}}{rT_{21}},\ r\in (0,1)$.
We have $rT_{21}=V(r)D_{rT_{11}}$,
where
      \begin{equation} \label{4.44}
V(r) :=rVD_{T_{11}}D_{rT_{11}}^{-1}=rV(I-T_{11}^*T_{11})^{1/2}\cdot
(I-r^2T_{11}^*T_{11})^{-1/2}.
      \end{equation}
Let us define the operator $U(r)$ by  \eqref{4.21} with  $V$
replaced by $V(r)$, but not replacing $M$. Then the operator
      \begin{equation} \label{4.45}
Q_0(r) :=D_{V^*(r)}^{-1}\cdot \bigl(I-V(r)U(r)\bigr)\cdot D_{U(r)}^{-1}
       \end{equation}
is bounded and $Q_0(r)\cos\ff$ is contractive.
We set $Q(r):=\overline{Q_0(r)}$ and note that
$Q(r)\cos\ff \in C(\pi/2)$.

Next, starting with $U(r)$ we define a dual pair of
Hermitian contractions
$\{rT_1,T_2(r)\}$ by  setting
$$
T_{21}(r):=D_{rT_{11}}U(r) \qquad \text{and}\qquad
T_2(r):=\binom{rT_{11}}{T_{21}^*(r)}.
$$
Denote by $T_K(r) (\in\Ext_{\{rT_1,T_2(r)\}})$ the
extension  of $\{T_1,T_2\}$ defined by the same operator $K$,
as the extension $T_K\in \Ext_{\{T_1,T_2\}},$  that is
     \begin{equation} \label{4.46}
T_K(r) :=
\begin{pmatrix}
rT_{11}&T_{12}(r)\\
rT_{21}&T_{22}(r)\end{pmatrix},\qquad
T_{22}(r) :=-rV(r)T_{11}U(r)+D_{V^*(r)}KD_{U(r)}.
    \end{equation}
Since $0\in\rho(D_{rT_{11}})\cap \rho(D_{V(r)})$ and
$Q(r)\cos\ff \in C(\pi/2)$, then for the operator-matrix
      \begin{equation}\label{4.46A}
(G_{ij}^{\pm}(r))_{i,j=1}^2 :=I-T_K^*(r)T_K(r)\pm i\ctg\ff(T_K(r)-T_K^*(r))
       \end{equation}
 equality \eqref{4.32} is already proved in the previous step,
that is
        \begin{equation} \label{4.48}
    \begin{split}
&\sin^2\ff\cdot [G_{22}^{\pm}(r)-(G_{12}^{\pm}(r))^*G_{11}^{-1}(r)G_{12}(r)]\\
&=D_{U(r)}\cdot [I-(K^*\sin\ff\pm iQ^*(r)\cos\ff)
\cdot (K\sin\ff\mp iQ(r)\cos\ff)]\cdot D_{U(r)}.
\end{split}
       \end{equation}

It remains to justify the possibility to pass to the limit
in \eqref{4.48} as $r\to 1$.
We may assume without rstriction of generality that $\ker G_{11}=\{0\}$.
Then, as it follows from
\eqref{4.34}, \eqref{4.35} and \eqref{4.38},
      \begin{equation} \label{4.49}
G_{11}^{1/2}=U_1D_VD_{T_{11}}=D_{T_{11}}D_VU_1^*,\quad
G_{11}^{-1/2}G_{12}^{\pm}=-U_1X_{\pm}^*,\quad
\bigl(G_{11}^{-1/2}G_{12}^{\pm}\bigr)^* = -X_{\pm}U_1^*,
      \end{equation}
where the operator $U_1$ is unitary.

Further, introducing the operators
      \begin{equation}\label{4.49A}
Y(r) :=(\sin^2\ff\cdot D_{V(r)}^2 +\ cos^2\ff)^{-1/2} \quad
\text{and}\quad  Y_*(r) :=(\sin^2\ff\cdot D_{V^*(r)} +\cos^2\ff)^{-1/2},
      \end{equation}
one derives from the definition of the operator $U^*(r)\ (r<1)$ that
      \begin{equation} \label{4.50}
U^*(r)-V(r)=\sin\ff D_{V^*(r)}Y_*(r)\cdot
\bigl(M^*-r\sin\ff\cdot V(r)\bigr)\cdot Y(r)D_{V(r)}.
      \end{equation}
Next, we define  $G_{11}^{1/2}(r)$ and  $G_{21}^{\pm}(r)$
by \eqref{4.34}  and \eqref{4.35} with $V(r), \ U(r)$ and $rT_{11}$
in place of $V, \ U$ and $T_{11}$ respectively.
Further, similarly to definition \eqref{4.38} of $X_{\pm}$ we set
        \begin{equation}\label{4.50A}
  X_{\pm}(r) := rU^*(r)T_{11}D_{V(r)}+D_{U(r)}K^*V(r)
 \mp i\cos\ff\cdot Y_*(r)D_{V*(r)} \bigl(M^*-r\sin\ff\cdot V(r)\bigr)Y(r).
      \end{equation}
Combining  these  definitions we arrive at the relations
      \begin{equation} \label{4.51}
 G_{11}^{1/2}(r)= U_1(r)D_{V(r)}D_{rT_{11}}  \quad \text{and}\quad
G_{21}^{\pm}(r)=X_{\pm}(r)D_{V(r)}D_{rT_{11}},
      \end{equation}
which are analogous to that of  \eqref{4.49}.
Here  $U_1(r),\ r\in (0,1),$ is a family of unitary operators.
Hence       
     \begin{equation} \label{4.52}
-(G_{12}^{\pm}(r))^*G_{11}^{-1/2}(r)=X_{\pm}(r)U_1^*(r),\qquad
-G_{11}^{-1/2}G_{12}^{\pm}(r)=U_1(r)X_{\pm}(r).
     \end{equation}
It follows from  \eqref{4.44}  that
$s-\lim_{r\to 1}V(r)=V\ $  and  $s-\lim_{r\to1}V^*(r)=V^*.$
Hence and taking into account \eqref{4.49A} we get
     \begin{eqnarray} \label{4.53}
  \begin{aligned}
    s-\lim_{r\to 1}D_{V(r)} &= D_V, &\qquad
    s-\lim_{r\to 1}D_{V^*(r)} &= D_{V^*}, \\
    s-\lim_{r\to 1}Y(r) &= Y, &
    s-\lim_{r\to 1}Y_*(r) &= Y_*.
  \end{aligned}
      \end{eqnarray}
Relations \eqref{4.50}, \eqref{4.53} and \eqref{4.37} yield
     \begin{equation}\label{4.53A}
 s-\lim_{r\to 1}U(r)=U,\qquad
 s-\lim_{r\to 1}U^*(r)=U^*,\qquad
s-\lim_{r\to 1}D_{U(r)}=D_U.
     \end{equation}
It follows from \eqref{4.46} and \eqref{4.46A} that
      \begin{equation}\label{4.53B}
s-\lim_{r\to 1}G_{22}(r)=G_{22}.
        \end{equation}
Further, \eqref{4.50A} and \eqref{4.38} yield
$ s-\lim_{r\to 1}X_{\pm}(r)=X_{\pm}\  $ and
$\  s-\lim_{r\to 1}X_{\pm}^*(r)=X_{\pm}^*.$
Therefore combining relations \eqref{4.48} with \eqref{4.52}
and taking into account the obvious identities
$U_1^*(r)U_1(r) = U_1^*U_1 =I$ we arrive at
       \begin{equation}\label{4.53C}
  s-\lim_{r\to 1}(G_{12}^{\pm}(r))^*G_{11}^{-1}(r)G_{12}^{\pm}(r)
    = s-\lim_{r\to 1}X_{\pm}(r)X_{\pm}^*(r)
    = X_{\pm}X_{\pm}^*
    = (G_{11}^{-1/2}G_{12}^{\pm})^*G_{11}^{-1/2}G_{12}^{\pm}.
        \end{equation}
Relations \eqref{4.53B} and \eqref{4.53C} allow us to pass to the
limit in left-hand side of \eqref{4.48} as $r\to 1$. So, it remains
to justify passage to the limit in the right hand side of \eqref{4.48}.
In turn it suffices to prove the relations
    \begin{equation} \label{4.54}
s-\lim_{r\to 1}Q(r)D_{U(r)}=QD_U\qquad \text{and}\qquad
 s-\lim_{r\to 1}D_{U(r)}Q^*(r)=D_UQ^*.
     \end{equation}
We derive from \eqref{4.50} and \eqref{4.45} that
      \begin{eqnarray}\label{4.54a}
Q(r)D_{U(r)} &=& D_{V^*(r)}^{-1}\bigl(I-V(r)U(r)\bigr) \\
\nonumber
&=& \{I-\sin\ff Y(r)V(r)
\bigl(M-r\sin\ff\cdot V^*(r)\bigr)Y_*(r)\}D_{V^*(r)}.
         \end{eqnarray}
It follows from \eqref{4.53} that there exists the limit of
the right-hand side  of \eqref{4.54a} as $r\to 1.$
Hence there exist the limit of the left-hand side
of \eqref{4.54a} as $r\to 1.$
Moreover, the first of relations
\eqref{4.54} is now implied by  \eqref{4.54a} and similar formula
for  $QD_U$ which follows from \eqref{4.37}.
The second formula in \eqref{4.54}
may be proved similarly.

Finally, passing to the limit in \eqref{4.48} as $r\to 1$
and taking into account \eqref{4.53B}, \eqref{4.53C} and
\eqref{4.54} we arrive at  \eqref{4.32}.
Relation \eqref{4.31} may  be proved in just the same way.
      \end{proof}

{\bf 4.2. Descriptions of the classes
$C(\pi/2;\kk_{\pm})$  and $T\in C(\pi/2;\gS^{\pm})$}

Here we present  some corollaries from Theorem \ref{th4.3}. To
formulate them we need some definitions and an elementary lemma.

Let $\kk_-(\gt)$ be the number of negative squares of the symmetris quadratic
form ${\gt}$, that is the maximum dimensions of the "negative"
linear manifolds
    \begin{equation*}
L\_=\{f\in\cD(\gt)\setminus \{0\}:{\gt}[f]<0\}\cup \{0\}.
    \end{equation*}
For any  selfadjoint operator $T=T^*\in\cC(\gH)$
with the resolution of the
identity $E_T(\cdot)$ we let
$T_-:=E_T(-\infty, 0)T$ and
$\kk\_(T):=\text{dim}(\ran T_-)= \text{dim}E_T(-\infty,0)\gH$.
If the form ${\gt}$ is closed and $T$ is  the operator
associated with it, ${\frak t}={\frak t}_T,$
(see \cite{Ka}) then by virtue of the minimax principle
$\kk\_(\gt)=\kk\_(T)$.

Next we define the classes
$C(\pi/2;\kk_{\pm})$  and $T\in C(\pi/2;\gS^{\pm})$.
      \begin{definition}\label{def2.5}
Let $\kk\in{\Bbb Z}_+$ and $\gS$ be a
two-sided ideal in $[\gH]$. We write

(a) $T\in C(\pi/2;\kk)$ if $T\in [\gH,\gH']$ and $\kk\_(I-T^*T)=\kk;$

(b) $T\in C(\pi/2;\gS)$ if $T\in [\gH,\gH']$ and $(I-T^*T)\_\in\gS$.
    \end{definition}
          \begin{definition}\label{def3.4}
Let $\ff\in[0,\pi/2]$, $\kk^{\pm}\in{\Bbb Z}^+,$ and let $\gS^{\pm}$
be two-sided ideals in $[\gH]$.  An operator $T(\in [\gH])$ is put

  (a)
  in the class $C_{\gH}(\ff;\kk^{\pm})$ with $\ff \in (0,\pi/2],$ if
  \begin{equation}\label{3.18}
    T\sin\ff\pm i\cos \ff\cdot I\in C_{\gH}(\pi/2;\kk^{\pm});
  \end{equation}

 (b)
in the class $C_{\gH}(\ff;\gS^{\pm})$  with $\ff \in (0,\pi/2],$ if
$$
T\sin\ff\pm i\cos\ff\cdot I\in C_{\gH}(\pi/2;\gS^{\pm}).
$$

(c)
in the class  $C_{\gH}(0;\kk)\  $ ($C_{\gH}(0;\gS)$), if $T=T^*$ and
$\kk\_(T)=\kk\quad (T_-\in \gS)$.

We write $C_{\gH}(\ff;\kk)$ and $C_{\gH}(\ff;\gS)$ in place of
$C_{\gH}(\ff;\kk^{\pm})$ and $C_{\gH}(\ff;\gS^{\pm})$ respectively
if $\kk:=\kk^+=\kk^-$ and $\gS:=\gS^+=\gS^-$;
    \end{definition}
Observe that the class $C_{\gH}(\pi/2;\kk^{\pm})$
is not empty only if $\kk^+=\kk^-$.
Some properties  of the class $C(\pi/2;\gS_{\infty})$
can be found in \cite{M4}.

     \begin{lemma}\label{lem2.5}  \cite{M3,HMS}
Let $T_1= \binom{T_{11}}{T_{21}}(\in[\gH_1,\gH])$ be
a nonnegative symmetric operator $(\equ T_{11}\ge 0)$
admitting a bounded nonnegative selfadjoint extension
and let $T(\in [\gH])$ be any selfadjoint extension
of $T_1$ with the block-matrix
representation  $T=T^*=(T_{ij})_{i,j=1}^2$ with respect to the orthogonal
decomposition $\gH=\gH_1\oplus\gH_2.$  Then

(i) $\Re(T_{11}^{1/2})\supset \Re(T_{12})$ and
the operator $S:=T_{11}^{-1/2}T_{12}$
is well-defined and bounded;

(ii) $\kk\_(T)=\kk\_(T_{22}-S^*S).$
In particular, $T\ge0$ iff $T_{22}-S^*S\ge 0.$
      \end{lemma}

Now we are ready to present the corollaries.
      \begin{corollary} \label{cor4.1C}
Let $\dup$ be a dual pair of symmetric  contractions,
$T_K\in \Ext_{\{T_1,T_2\}},\ $   $\ff\in[\ff_1,\pi/2]$ and $\ff_1>0$.
Then the following equivalences
are valid:
    \begin{equation}\label{4.54kappa}
T_K\in C_{\gH}(\ff;\kk^{\pm})\equ K\sin\ff\mp iQ\cos\ff\in C(\pi/2;\kk^{\pm}).
    \end{equation}
        \end{corollary}
    \begin{proof}
Let as in Theorem \ref{th4.3}
$$
(G_{ij}^{\pm})_{i,j=1}^2:=G^{\pm}:=I-T_K^*T_K\pm i\ctg \ff(T_K-T_K^*).
$$
The operators
$G_0^{\pm} := \binom{G_{11}}{G_{21}^{\pm}}$
are nonnegative $(\equ G_{11}\ge 0).$
Moreover, both of them  admit  bounded nonnegative
selfadjoint extensions. For example, the operator
    \begin{equation*}
\begin{pmatrix}
G_{11}^{1/2}&0\\B_+^*&0\end{pmatrix}\,
\begin{pmatrix}
G_{11}^{1/2}&B_+\\0&0\end{pmatrix}=
\begin{pmatrix}
G_{11}&G_{12}^+\\G_{21}^+&B_+^*B_+
\end{pmatrix}\ge 0
    \end{equation*}
with a bounded   $B_+=G_{11}^{-1/2}G_{12}^+$ is
a nonnegative extension of $G_0^+$.
Therefore combining
Lemma \ref{lem2.5} with Theorem \ref{th4.3}
(see equality \eqref{4.32}) we get
$$
\kk\_(I-T_K^*T_K)=\kk\_(G_{22}^{\pm}-(G_{11}^{-1/2}G_{12}^{\pm})^*
(G_{11}^{-1/2}G_{12}^{\pm}))=\kk\_(I-K^*_{\mp}K_{\mp}),
$$
where $K_{\pm}:=K sin\ff\pm iQ cos\ff$.
     \end{proof}
       \begin{corollary}\label{cor4.2C}
Let $\dup$ be a dual pair of symmetric  contractions
and $0\in \rho(D_{T_{11}})\cap \rho(D_V).$
Suppose additionally that $\ff\in[\ff_1,\pi/2]$ and $\ff_1>0$.
Then the following implications hold
     \begin{equation}\label{4.54ideal}
K\sin\ff\pm Q\cos\ff \in C(\pi/2;\gS^{\pm}) \then
T_K\in C_{\gH}(\ff;\gS^{\pm}).
     \end{equation}
If additionally $0\in \rho(D_U)$ then implications
\eqref{4.54ideal} turns into the equivalences.
        \end{corollary}
        \begin{proof}
The required assertion immediately follows from \eqref{4.32}
and the identity
     \begin{equation*}
     \begin{pmatrix}
I&0\\-G_{21}^{\pm}G_{11}^{-1}&I\end{pmatrix}\,
\begin{pmatrix} G_{11}&G_{12}^{\pm}\\
G_{21}^{\pm}&G_{22}^{\pm}\end{pmatrix}\,
\begin{pmatrix} I&-G_{11}^{-1} G_{12}^{\pm}\\
0&I\end{pmatrix}=
\begin{pmatrix} G_{11}&0\\
0&G_{22}^{\pm}-G_{21}^{\pm}G_{11}^{-1}G_{12}^{\pm}
   \end{pmatrix}.
     \end{equation*}
\end{proof}
     \begin{remark}\label{rem2.2A}
(i) Let $\ff=\pi/2.$ Then both relations \eqref{4.54kappa}
and \eqref{4.54ideal} are simplified and take the form
      \begin{equation}\label{2.51}
K\in C(\pi/2;\kk)\equ T_K\in C(\pi/2;\kk),
       \end{equation}
     \begin{equation}\label{2.52}
(I-K^*K)\_\in\gS \then (I-T^*_KT_K)\_ \in\gS.
      \end{equation}
Both relations have been established in \cite{M2, M4}
for any (not necessary symmetric) dual pair of contractions.

(ii) Let $\ff_1=\arccos(\|Q_0\|^{-1})= 0.$
Then $\|Q_0\|=1$ and by Remark \ref{rem4.3}
(see \eqref{4.19A}) $U=V^*$, that is $Q=I$ and $T_1=T_2$.
In this case description of the sets
$\Ext_{T_1}(\ff;\kk)= \Ext_{\{T_1,T_1\}}(\ff;\kk)$ and
$\Ext_{T_1}(\ff;\gS)=\Ext_{\{T_1,T_1\}}(\ff;\gS), \quad
\ff\in[0,\pi/2]$,
can  easily be derived from Corollaries \ref{cor4.1C}
and \ref{cor4.2C}.    Now in place of relations
\eqref{4.54kappa} and \eqref{4.54ideal} we have
    \begin{equation}\label{2.51A}
K \in C_{\cH}(\ff;\kk^{\pm}) \equ  T_K\in C_{\gH}(\ff;\kk^{\pm}),
    \end{equation}
       \begin{equation}\label{2.52A}
K \in C_{\cH}(\ff;\gS^{\pm}) \then  T_K \in C_{\gH}(\ff;\gS^{\pm}),
        \end{equation}
where $\cH={\overline \ran}D_U = {\overline \ran}D_{V^*}.$
Both formulas  have earlier  been obtained in \cite{M2, M4}.
Note also, that if $\kk^{\pm}=0$ then formula
\eqref{2.51A}  gives one more proof of Corollary
\ref{cor4.1}.
        \end{remark}
      \begin{definition}\label{def3.3}
Let $\ff\in[0,\pi/2),\  \kk^+\in{\Bbb Z}_+,\  \gS^{\pm}$ two-sided
ideals in $[\gH],$ and  $B\in \cC(\gH).$
Let further, the quadratic forms
     \begin{equation*}
\gt_{\pm}[f]=\re (Bf,f)\pm\ctg\ff\cdot\im(Bf,f),\quad f\in \dom B,
    \end{equation*}
be semibounded below
and  $B^{\pm}$  the linear operators
associated with their closures (the closability
of the form $\gt_{\pm}$ is a consequence of their
semiboundness (see \cite{Ka})).
We write

(a) $B\in S_{\gH}(\ff;\kk^{\pm})$, if $\rho(B)\not=\eset$ and
$\kk(B^{\pm})=\kk^{\pm}$;

 (b) $B\in S_{\gH}(\ff;\gS^{\pm})$, if $(B^{\pm})\_\in\gS^{\pm}$ and
$\rho(B)\not=\eset$.

A closed linear relation $\theta$ in $\gH$ is also put in the class
$S_{\gH}(\ff;\kk^{\pm})\  (S_{\gH}(\ff;\gS^{\pm}))$
if Re$(f',f)\ge \beta\|f\|^2$\  for all $\{f,f'\}\in \theta$
(with some $\beta\in {\Bbb R}$) and
its operator part is in $S_{\gH}(\ff;\kk^{\pm})\
(S_{\gH}(\ff;\gS^{\pm}))$.

It is clear that $S_{\gH}(\ff;0)$ coinsides with $S_{\gH}(\ff)$.
      \end{definition}
It is clear that the classes $C_{\gH}(\ff;\kk^{\pm})$  and
$S_{\gH}(\ff;\kk^{\pm})$  are connected by means of the linear
fractional transformation  \eqref{3.2}.
  The same is also true for the classes $C_{\gH}(\ff;\gS^{\pm})$ and
  $S_{\gH}(\ff;\gS^{\pm})$.


{\bf 4.3. Shorted operators.}
Here we present  two  additional corollaries
from Theorem \ref{th4.3} complementing Theorem \ref{4.1}.
For this purpose we recall some well-known results
and the definition of a  shorted operator.
        \begin{definition}\label{def.2.6}(\cite{Kr1})
For any nonnegative operator $A(\in[\gH])$
and a subspace $\gN(\subset \gH)$ there exists
the largest element in the set
of all bounded operators not exceeding $A$ and
annihilating $\gN^{\perp}=\gH\ominus\gN$.
This element, is denoted by $A_\gN$ and is called the
shorted to $\gN$ operator.

The transformation $A\to A_\gN$ is called the Krein transformation.
   \end{definition}
   \begin{lemma}\label{lem4.1A} (\cite{Kr1, Sh, HMS}).
Let $A=(A_{ij})_{i,j=1}^2$ be a block-matrix
representation of an operator
$A\ge 0\,(A\in[\gH])$
with respect to the decomposition $\gH=\gH_1\oplus \gN$. Then the
shorted to $\gN$ operator $A_\gN$ is of the form
    \begin{equation}\label{2.54}
A_\gN=
\begin{pmatrix}
0&0\\
0&A_{22}-S^*S
\end{pmatrix},\qquad S=A_{11}^{-1/2}A_{12}.
     \end{equation}
If in addition $0\in \rho(A_{11})$ then $S^*S=A_{21}A_{11}^{-1}A_{12}$.
      \end{lemma}
      \begin{corollary}\label{cor4.5}(\cite{Kr1, KO, Sh}).
Let $\gH=\gH_1\oplus\gN,\ A\in[\gH]$
and $A\ge 0$. Then
    \begin{equation}\label{2.55}
\inf_{g\in\gH_1}(A(f-g),f-g)=(A_\gN f,f),\qquad \  f\in \gH.
    \end{equation}
       \end{corollary}
       \begin{corollary}\label{cor4.3}
Let $\dup$ be a dual pair of symmetric  contractions in
$\gH=\gH_1\oplus\gH_2,\ $ $\ff\in [\ff_1,\pi/2]$ where
$\ff_1:=\arccos(\|{Q_0}\|^{-1})$
and  $T_K\in \Ext_{\{T_1,T_2\}}(\ff)$.
Then the shorted to $\gN:=\gH_2$ operators
    \begin{equation}\label{2.55G}
G^{\pm}=I-T_K^*T_K\pm i\ctg \ff(T_K-T_K^*)\quad \text{and} \quad
S^{\pm}=I-T_KT_K^*\pm i\ctg\ff(T_K-T_K^*)
    \end{equation}
have the following form
    \begin{equation}\label{4.54A}
(G^{\pm})_{\gN}=
\begin{pmatrix}0&0\\
0&\sin^{-2}\ff\cdot D_U[I-(K^*\sin\ff\pm i Q^*\cos\ff)(K\sin\ff\mp iQ\cos\ff)]D_U
\end{pmatrix}
      \end{equation}
and
    \begin{equation}\label{4.54B}
(S^{\pm})_{\gN}=
\begin{pmatrix}
0&0\\
0&\sin^{-2}\ff\cdot D_{V^*}[I-(K\sin\ff\mp iQ\cos\ff)
(K^*\sin\ff\pm iQ^*\cos\ff)]D_{V^*}
\end{pmatrix}.
     \end{equation}
        \end{corollary}
     \begin{proof}
One deduces the proof combining Theorem \ref{th4.3} with
Lemma  \ref{lem4.1A}.
   \end{proof}
       \begin{corollary}\label{cor4.3A}
Let $\dup$ be a dual pair of Hermitian contractions in
$\gH=\gH_1\oplus\gH_2,\ $ and
$T_K\in \Ext_{\{T_1,T_2\}}(\ff)$.
Suppose additionally that
$\ff_1=\arccos(\|Q_0\|^{-1})<\pi/2,\ $
$\ff\in[\ff_1, \pi/2]$
and $\gS^{\pm}$ are two-sided ideals in $[\gH]$. Then

(i)  $\Ext_{\{T_1,T_2\}}(\ff)\not =\emptyset$ if and only if
$\ff\in [\ff_1, \pi/2];$

(ii) The following equivalence  holds
    \begin{equation}\label{4.54C}
T_K\in \Ext_{\{T_1,T_2\}}(\ff)\equ  K_{\pm}:=K\sin\ff\pm i Q\cos\ff\in C(\pi/2);
       \end{equation}

(iii)  The following implications  hold  with $\gN:=\gH_2$
      \begin{equation}\label{4.54D}
D^2_{K_{\mp}}\in \gS^{\pm}  \then (G^{\pm})_{\gN}\in \gS^{\pm}, \qquad
D^2_{K_{\mp}^*}\in \gS^{\pm}  \then (S^{\pm})_{\gN}\in \gS^{\pm}.
         \end{equation}
If additionally $0\in \rho(D_U)$ then implications
\eqref{4.54D} turn into the equiavalences.
         \end{corollary}
       \begin{proof}

(i)-(ii)
By definition $T_K\in C_{\gH}(\ff)$ if and only if
$G^{\pm}\in C_{\gH}(\pi/2),$ where $G^{\pm}$ are defined by
\eqref{2.55G}.
Note that
$G^{\pm}_{11}=G_{11}=I-T^*_{11}T_{11}-T^*_{21}T_{21}\ge 0$ since
$T_1$ is a contraction.
Therefore by Sylvester criterion (see Lemma \ref{lem2.5})
$G^{\pm}\ge 0$ iff $(G^{\pm})_{\gN}\ge 0$
with $\gN=\gH_2$.
Combining this inequality with \eqref{4.54A}
we arrive at  equivalence \eqref{4.54C}.

Hence, if $\Ext_{\{T_1,T_2\}}(\ff)\not =\emptyset$ then
$Q\cos\ff\in C(\pi/2)$, that is
$\ff\in[\ff_1,\pi/2]$.
Conversly, if $\ff\in [\ff_1,\pi/2]$ then
$Q\cos\ff\in C(\pi/2)$ and
the operator $T_0$, that is $T_K$ with $K=0$, belongs to
$\Ext_{\{T_1,T_2\}}(\ff)$.

(iii) This statement is immediately implied by formulas
\eqref{4.54A} and \eqref{4.54B}.
         \end{proof}
   \begin{remark}\label{rem4.5A}
(i) Suppose that  in Corollary \ref{cor4.3} $\ff=\pi/2.$
Then $G^{\pm}=I-T_K^*T_K= D_{T_K}^2$ and
$S^{\pm}=I-T_KT_K^* =D_{T_K^*}^2.$  Now
formulas \eqref{4.54A} and \eqref{4.54B}
are simplified and take the form
    \begin{equation}\label{2.56}
(D_{T_K}^2)_{\gN}=
\begin{pmatrix} 0&0\\0&D_UD_K^2D_U
\end{pmatrix},\qquad
(D_{T_K^*}^2)_{\gN'}=
\begin{pmatrix}
0&0\\0&D_{V^*}D_K^2D_{V^*}
\end{pmatrix}.
    \end{equation}
where $\gN:=\gH_2$ and $\gN':=\gH'_2$.
Both formulas have earlier been obtained in \cite{M2, M4}
for any (not necessary symmetric) dual pair of contractions.

(ii)  Corollary \ref{cor4.3A} (iii) complements
Theorem \ref{th4.1}.
Moreover, Corollary \ref{cor4.3A} gives another proof
of Theorem  \ref{th4.1}.
Indeed, in the case
$0\in\rho(D_U)\cap\rho(D_{V^*})$ the proof of
Theorem \ref{th4.3} does not
depend on Theorems \ref{th4.1} and \ref{th4.2}.
The proof of equivalence \eqref{4.54C}
without the additional assumption  $0\in\rho(D_U)\cap\rho(D_{V^*})$
can easily be obtained
by considering the family $\{rT_1, rT_2\},\ r\in (0,1)$,
of dual pairs  of contractions and
passage to the limit as $r\to 1$
(cf. the proof of Theorem \ref{th4.2}).
     \end{remark}
\section{ Completions of a special triangular operator-matrix.}

{\bf 5.1. A complement to the S. Nagy and C. Foias result.}
Here we describe the operators $T_{12}$
completing the incomplete contractive operator block-matrix
    \begin{equation} \label{4.55}
\begin{pmatrix}
T_{11}&*\\
0&T_{22}
\end{pmatrix}
    \end{equation}
to form an operator matrix of some class.

We start with the following S. Nagy and C. Foias result.
       \begin{proposition}\label{pr4.1}\cite{SNF1}
Let
$\gH=\gH_1\oplus\gH_2=\gH_1'\oplus\gH_2'$
and let  $T_{jj} (\in [\gH_j,\gH_j'])$ be a contraction,
$ j\in\{1,2\}$.
Then the family of operators $T_{12}\in [\gH_2,\gH_1'],$
completing the block-matrix \eqref{4.55}
to a contractive matrix
$T=(T_{ij})_{i,j=1}^2\in[\gH]$,
forms an operator ball $B(0; D_{T_{11}^*},D_{T_{22}})$, that is
     \begin{equation*}
T=
\begin{pmatrix}
T_{11}&T_{12}\\
0&T_{22}
\end{pmatrix}
\in C_\gH(\pi/2)\equ T_{12}=D_{T_{11}^*}KD_{T_{22}},\quad
 \|K\|\le 1,\quad  K\in [\cH_2,\cH_1].
    \end{equation*}
         \end{proposition}
       \begin{proof}
Let at first $0\in \rho(D_{T_{22}})$.
Then according to the Sylvester criterion the equivalence
    \begin{equation} \label{4.56}
I-T^*T\ge 0\equ D_{T_{22}}^2-Z^*(I+T_{11}D_{T_{11}}^{-2}T_{11}^*)Z\ge 0,
     \end{equation}
holds true with $Z:=T_{12}$.
By Lemma \ref{lem2.2} the set of solutions of
\eqref{4.56} forms an operator ball.
Observing that $I+T_{11}D_{T_{11}}^{-2}T_{11}^2=
D_{T_{11}^*}^{-2}$ and applying Lemma \ref{lem2.2} to
\eqref{4.56} we arrive at the
required relation
   \begin{equation*}
T_{12}=Z=D_{T_{11}^*}KD_{T_{22}},\qquad \|K\|\le 1,\quad
 K\in [\cH_2,\cH_1].
    \end{equation*}
We may easily free ourselves of the condition $0\in\rho(D_{T_{22}})$
by virtue of passage to the limit.
        \end{proof}

Thus the contractive "completions" of the matrix \eqref{4.55}
are of the form
     \begin{equation} \label{4.57}
T_K=
\begin{pmatrix}
T_{11}&D_{T_{11}^*}KD_{T_{22}}\\
0&T_{22}
\end{pmatrix}
    \end{equation}
with $\|K\|\le 1,\ K\in [\cH_2.\cH_1]$.

Let us now consider "completions" of the incomplete
block-matrix \eqref{4.55}
of the form \eqref{4.57} not assuming the operator
$K$ to be a contraction.
    \begin{proposition}\label{pr4.2}
Let $\gH=\gH_1\oplus\gH_2=\gH'\oplus\gH_2',\
T_{jj}\in [\gH_j,\gH_j'],$ and  $ \|T_{jj}\|\le 1,\ j\in \{1,2\}.$
Assume that $T_K$ is an operator matrix of the
form \eqref{4.57}  with $K\in [\cH_2,\cH_1]$ and put
$G:=(G_{ij})_{i,j=1}^2=I-T_K^*T_K$
and $\ S:=(S_{ij})_{i,j=1}^2=I-T_KT_K^*.$
Then

(i) $\ran(G_{12})\subset \ran(G_{11}^{1/2})$  and
$\ran(S_{12})\subset \ran(S_{11}^{1/2}),$  hence
the operators $G_{11}^{-1/2}G_{12}$ and $S_{11}^{-1/2}S_{12}$
are well defined and bounded;

(ii) the following identities are valid:
    \begin{equation} \label{4.58}
G_{22}-(G_{11}^{-1/2}G_{12})^*(G_{11}^{-1/2}G_{12})=
D_{T_{22}}(I-K^*K)D_{T_{22}},
      \end{equation}
      \begin{equation}\label{4.58A}
S_{22}-(S_{11}^{-1/2}S_{12})^*(S_{11}^{-1/2}S_{12})=
D_{T_{11}^*}(I-KK^*)D_{T_{11}^*}.
    \end{equation}
         \end{proposition}
        \begin{proof}
Imposing the condition $0\in \rho(D_{T_{11}})$ we have
      \begin{eqnarray*}
  G_{22}-G_{12}^*G_{11}^{-1}G_{12} \\
    = D_{T_{22}}(I-K^*D_{T_{11}^*}K)D_{T_{22}}-D_{T_{22}}
      K^*D_{T_{11}^*}T_{11}D_{T_{11}}^{-2}T_{11}^*D_{T_{11}^*}KD_{T_{22}} \\
    = D_{T_{22}}(I-K^*K+K^*T_{11}T_{11}^*K-K^*T_{11}T_{11}^*K)D_{T_{22}}
    = D_{T_{22}}D_K^2D_{T_{22}}.
         \end{eqnarray*}
We may free ourselves of the condition
$0\in\rho(D_{T_{11}})$ by passage to the limit
just like it was done in the proof of
Theorem \ref{th4.2}.
Equality \eqref{4.58A} may be proved similarly.
      \end{proof}
      \begin{corollary}\label{cor4.4}
Suppose that conditions of Proposition \ref{pr4.2}
are satisfied. Then

(i) the following equivalence holds
    \begin{equation*}
K\in C(\pi/2;\kk)\cap [\cH_2,\cH_1]\equ  T_K\in C_{\gH}(\pi/2;\kk);
     \end{equation*}

(ii) if in addition $0\in\rho(G_{11})$, then
for any two-sided ideal $\gS$ in $[\gH]$
the following implication holds
      \begin{equation*}
K\in C(\pi/2;\gS)\then T_K\in C_{\gH}(\pi/2;\gS).
       \end{equation*}
This implication  turns into the equivalence
if additionally  $0\in\rho(D_{T_{22}}).$
      \end{corollary}
        \begin{corollary} \label{cor4.5}
 Let $T_K$ be a contraction of the form \eqref{4.57}.
Then the  operators $G:=D_{T_K}^2$
and $S:=D_{T_K^*}^2$ shorted to $\gN:=\gH_2$ and $\gN':=\gH_2'$
respectively have the form
     \begin{equation*}
G_{\gN}=
\begin{pmatrix}
0&0\\
0&D_{T_{22}}D_K^2D_{T_{22}}
\end{pmatrix}
\qquad \text{and} \qquad
S_{\gN'}=\begin{pmatrix}
0&0\\
0&D_{T_{11}^*}D_{K^*}^2D_{T_{11}^*}
\end{pmatrix}.
     \end{equation*}
\end{corollary}

Corollaries \ref{cor4.4}  and \ref{cor4.5}  may be derived
from Proposition \ref{pr4.2} just like
Corollaries \ref{cor4.1C}, \ref{cor4.2C} and
\ref{cor4.3} from  Theorem \ref{th4.3}.

{\bf 5.2. A solution to   Yu.~L.~Shmul'yan's problem.}

The following  proposition provides an answer
to  the Yu.~L.~Shmul'yan  question.
       \begin{proposition}\label{pr4.3}
Let $\gH=\gH_1\oplus\gH_2,\ T_{jj}\in C_{\gH_j}(\ff),\
j\in \{1,2\},\ $ and   $\ff\in (0,\pi/2)$. Then

(i) there exist  contractions
$U_{\ff}=U_{\ff}^*\in C_{\gH_1}(\pi/2)$ and
$V_{\ff}=V_{\ff}^*\in C_{\gH_2}(\pi/2)$
such that
     \begin{equation} \label{4.60}
2\ctg\ff(\im T_{11})=D_{T_{11}^*}U_{\ff}D_{T_{11}^*}\quad \text{and}
\quad 2\ctg\ff(\im T_{22})
=D_{T_{22}}V_{\ff}D_{T_{22}}\ ;
     \end{equation}

(ii) the following equivalence  holds  true
     \begin{eqnarray} \label{4.61}
T=
\begin{pmatrix}
T_{11}&T_{12}\\
0&T_{22}
\end{pmatrix}
\in C_{\gH}(\ff) \equ  \\
\nonumber
T_{12}=\sin\ff D_{T_{11}^*}KD_{T_{22}},\quad
(I\pm U_{\ff})^{-1/2}K(I\pm V_\ff)^{-1/2}\in C(\pi/2),\quad
K\in [\cH_2,\cH_1].
   \end{eqnarray}
      \end{proposition}
\begin{proof}
(i) Equalities \eqref{4.60}  have already been proved in
Proposition \ref{prop4extreme}.

(ii) The inclusion  $T\in C_\gH(\ff)$ means that
$T\sin\ff\pm i\cos\ff\cdot I\in C_\gH(\pi/2)$.
Since $T_{jj}\sin\ff\pm i\cos\ff\cdot I\in C_{\gH_j}(\pi/2),\
j\in \{1,2\},$  then by
Proposition \ref{pr4.1} the equivalences
     \begin{equation} \label{4.62}
T\sin\ff\pm i\cos\ff\cdot I\in C(\pi/2)\equ T_{12}\sin\ff=
R_l^{\pm}K_{\pm}R_r^{\pm},\quad  \|K_{\pm}\|\le 1,
    \end{equation}
hold true. Here $K_{\pm}\in [\cH_r^{\pm},\cH_l^{\pm}],\ \cH_r^{\pm}=
{\overline \ran}(R_r^{\pm}),\ \cH_l^{\pm}=
{\overline \ran}(R_l^{\pm})$,
and the operators $R_l^{\pm}$ and $R_r^{\pm}$
are defined by
     \begin{equation} \label{4.63}
\begin{split}
(R_l^{\pm})^2=\sin^2\ff D_{T_{11}^*}^2\pm \sin2\ff\cdot
(\text{Im}T_{11})=\sin^2\ff\cdot D_{T_{11}^*}(I\pm U_\ff)D_{T_{11}^*},
\\
(R_r^{\pm})^2=\sin^2\ff D_{T_{22}}^2\pm \sin2\ff\cdot
(\text{Im}T_{22})=\sin^2\ff\cdot D_{T_{22}}(I\pm V_\ff)D_{T_{22}}.
\end{split}
    \end{equation}
It is clear that
    \begin{equation} \label{4.64}
R_l^2:=(R_l^+)^2+(R_l^-)^2=2\sin^2\ff D_{T_{11}^*},\quad
R_r^2:=(R_r^+)^2+(R_r^-)^2=2\sin^2\ff D_{T_{22}}^2.
    \end{equation}
Further,  relations \eqref{4.63} yield polar
representations for the operators
$\sin\ff\cdot(I\pm U_\ff)^{1/2}D_{T_{11}^*}$
and $\sin\ff\cdot(I\pm V_\ff)^{1/2}D_{T_{22}}.$ Namely, we have
     \begin{equation} \label{4.65}
\sin\ff\cdot(I\pm U_\ff)^{1/2}D_{T_{11}^*}=U_{\pm}R_l^{\pm}\quad
\text{and} \quad
\sin\ff\cdot(I\pm V_\ff)^{1/2}D_{T_{22}}=V_{\pm}R_r^{\pm},
     \end{equation}
where $U_{\pm}$ and $V_{\pm}$ are partial isometries with initial spaces
$\cH_l^{\pm}$ and $\cH_r^{\pm}$ respectively.
We deduce the following equalities from \eqref{4.64} and \eqref{4.65}:
    \begin{equation} \label{4.66}
R_l^{\pm}= \frac1{\sqrt{2}}R_l(I\pm U_\ff)^{1/2}U_{\pm},
\qquad R_r^{\pm}=\frac1{\sqrt{2}}V_{\pm}^*(I\pm V_\ff)^{1/2}R_r.
    \end{equation}
Taking \eqref{4.66} into account we rewrite expression \eqref{4.62}
for $T_{12}\sin\ff$ in the form
   \begin{equation} \label{4.67}
T_{12}\sin\ff=
\frac12 R_l(I\pm U_\ff)^{1/2}U_{\pm}K_{\pm}V_{\pm}^*(I\pm V_\ff)^{1/2}R_r.
   \end{equation}
It follows from \eqref{4.67} that
    \begin{equation} \label{4.68}
(I+U_\ff)^{1/2}U_+K_+V_+^*(I+V_\ff)^{1/2}=
(I-U_\ff)^{1/2}U_-K_-V_-^*(I-V_\ff)^{1/2}.
    \end{equation}
Denoting  the operator in the left-hand side of
\eqref{4.68} by $K$ and taking into account \eqref{4.64} we arrive
at the following formula for $T_{12}:$
    \begin{equation*}
T_{12}=\frac1{2\sin\ff}R_lKR_r=\sin\ff D_{T_{11}^*}KD_{T_{22}},
    \end{equation*}
with $(I\pm U_\ff)^{-1/2}K(I\pm V_\ff)^{-1/2}\in C(\pi/2)$.
Here we have made use of the obvious  equivalences
$U_{\pm}K_{\pm}V_{\pm}^*\in C(\pi/2)\equ K_{\pm}\in C(\pi/2)$.

\end{proof}




\end{document}